\input ./style/arxiv-general.cfg
\documentclass[aap,MSNbibl,seceqn,dvips]{arximspdf}
\makeatletter
   \@ifpackageloaded{graphicx}{}{\usepackage{graphicx}}
\makeatother
\doi{10.1214/14-AAP1084}% Updated by VTEXPTS2LaTeX.exe, 22.12.2014
\volume{26}
\issue{1}
\pubyear{2016}
\firstpage{1}
\lastpage{44}
\docsubty{FLA}

\makeatletter
\newcommand{\rrvert}{\vert}
\newcommand{\rrVert}{\Vert}
\newcommand{\llvert}{\vert}
\newcommand{\llVert}{\Vert}
\renewcommand{\mid}{|}
\newtheorem{Theorem}{Theorem}[section]
\newtheorem{Proposition}[Theorem]{Proposition}
\newtheorem{Lemma}[Theorem]{Lemma}
\newproclaim{Remark}[Theorem]{Remark}
\newproclaim{ass}{Assumption}
\def\E{\mathbb{E}}
\def\tbf{\mathbf{t}}
\def\Dbf{\bolds{\Delta}}
\newcommand{\PP}{\mathbb{P}}
\newcommand{\RR}{{\mathbb R}}
\newcommand{\ub}{\underline{b}}
\makeatother

\begin{document}
\begin{frontmatter}

\title{The maximum maximum of a martingale with\\ given \textit{\lowercase{n}} marginals}
\runtitle{The maximum maximum of a martingale}

\begin{aug}
% Corresponding author: Jan Obloj - jan.obloj@maths.ox.ac.uk% Updated by VTEXPTS2LaTeX.exe, 29.12.2014 13:02
%by VTEXPTS2LaTeX.exe, 22.12.2014 10:37
\author[A]{\fnms{Pierre}~\snm{Henry-Labord\`ere}\ead[label=e1]{pierre.henry-labordere@sgcib.com}},
\author[B]{\fnms{Jan}~\snm{Ob\l{}\'{o}j}\corref{}\ead[label=e2]{Jan.Obloj@maths.ox.ac.uk}\thanksref{T1}},
\author[B]{\fnms{Peter}~\snm{Spoida}\ead[label=e3]{Peter.Spoida@maths.ox.ac.uk}\thanksref{T2}}
\and
\author[C]{\fnms{Nizar}~\snm{Touzi}\ead[label=e4]{nizar.touzi@polytechnique.edu}\thanksref{T3}}
\runauthor{Henry-Labord\`ere, Ob\l\'oj, Spoida and Touzi}
\affiliation{Soci\'et\'e G\'en\'erale, University of Oxford,\\
University of Oxford and Ecole Polytechnique Paris}
%\dedicated{}
\address[A]{P. Henry-Labord\`ere\\
Soci\'et\'e G\'en\'erale\\
Global Markets Quantitative Research\\
17 cours Valmy\\
92987 Paris-la D\'efense Cedex\\
France\\
\printead{e1}}
\address[B]{J. Ob\l{}\'{o}j\\
P. Spoida\\
Mathematical Institute\\
University of Oxford\\
ROQ, Woodstock Road\\
Oxford OX2 6GG\\
United Kingdom\\
\printead{e2}\\
\phantom{E-mail: }\printead*{e3}}
\address[C]{N. Touzi\\
Centre de Math\'ematiques Appliqu\'ees\\
Ecole Polytechnique Paris\\
91128 Palaiseau Cedex\\
France\\
\printead{e4}}
\end{aug}
\thankstext{T1}{Supported by ERC Starting Grant \textsc{RobustFinMath} 335421,
the Oxford-Man Institute of Quantitative Finance and St John's College
in Oxford.}
\thankstext{T2}{Supported by the Oxford-Man Institute of
Quantitative Finance and the DAAD.}
\thankstext{T3}{Supported by ERC Advanced Grant 321111 ROFIRM,
the Chair \textit{Financial Risks} of
the \textit{Risk Foundation} sponsored by Soci\'et\'e G\'en\'erale and the
Chair \textit{Finance and Sustainable Development} sponsored by EDF and CA-CIB.}

% HISTORY:
%
\received{\smonth{11} \syear{2013}}% Updated by VTEXPTS2LaTeX.exe,
%22.12.2014 10:37
%
\revised{\smonth{9} \syear{2014}}% Updated by VTEXPTS2LaTeX.exe,
%22.12.2014 10:37

% ABSTRACT
%
\begin{abstract}
We obtain bounds on the distribution of the maximum of a
martingale with fixed
marginals at finitely many intermediate times. The bounds are sharp and
attained by a solution
to $n$-marginal Skorokhod embedding problem in Ob\l\'oj and Spoida
[An iterated Az\'ema-Yor type embedding for finitely many marginals (2013) Preprint]. It follows that
their embedding maximizes the maximum among all other embeddings. Our
motivating problem is superhedging lookback options under volatility
uncertainty for an investor allowed to dynamically trade the underlying
asset and statically trade European call options for all possible
strikes and finitely-many maturities. We derive a pathwise inequality
which induces the cheapest superhedging value, which extends the
two-marginals pathwise inequality of Brown, Hobson and Rogers
[\textit{Probab. Theory Related Fields} \textbf{119} (2001) 558--578]. This
inequality, proved by elementary arguments, is derived by following
the stochastic control approach of Galichon, Henry-Labord\`ere and
Touzi [\textit{Ann. Appl. Probab.} \textbf{24} (2014) 312--336].
\end{abstract}

% KEYWORDS
\end{frontmatter}

\setcounter{footnote}{3}

%s1 #&#
\section{Introduction}

\subsection*{Probabilistic perspective}
The problem of controlling the maximum of a continuous martingale using
its terminal distribution has a long and rich history, starting with
Doob's maximal inequalities. In seminal contributions, Blackwell and
Dubins \cite{BD63}, Dubins and Gilat \cite{MR0494473} and Az\'ema and
Yor \cite{AzemaYor1,AzemaYor2} established that the distribution of
the maximum $X^*_T:=\sup_{t\leq T} X_t$ of a\vspace*{1pt} martingale $(X_t)$ is
bounded above, in stochastic order, by the so-called Hardy--Littlewood
transform of the distribution of $X_T$, and the bound is attained. This
led to series of studies on the possible distributions of $(X_T,X^*_T)$
including Gilat and Meilijson \cite{GM88}, Kertz and R\"osler \cite
{KR90,KertzRosler92b,KertzRosler93}, Rogers \cite{MR1217446} and
Vallois \cite{MR1268251}; see also Carraro, El Karoui and Ob\l\'oj
\cite{CKO12}.
More recently, such problems appeared very naturally within the field
of mathematical finance, as we explain below, which motivated further
developments. The original result was generalized to the case of a
nontrivial starting law in Hobson \cite{Hobsonmax} and to the case of
a fixed intermediate law in Brown, Hobson and Rogers \cite{brownhobsonrogers}.

In this paper, we generalize the above studies by controlling the
maximum using several intermediate marginals. We consider the case when
distributions of $X_{t_1},\ldots, X_{t_{n-1}}, X_{t_n}$ are given and
establish an upper bound on the distribution of the maximum
$X^*_{t_n}$. Our motivation comes from the mathematical finance problem
of robust superhedging of a lookback option. We apply a
general duality result from Possama\"i et al. \cite{Possamai-Royer-Touzi}
which converts the original problem into a min-max
calculus of variations problem where the Lagrange multipliers encode
the intermediate marginal constraints. The multipliers have in fact an
important financial interpretation as the optimal static positions in
Vanilla options, which reduce the risk induced by the derivative
security. Following Galichon, Henry-Labord\`ere and Touzi \cite{ght},
we apply stochastic control methods to solve the new problem
explicitly. The first step of our solution recovers the extended
optimal properties of the Az\'ema--Yor solution to the Skorokhod
embedding problem (SEP) obtained by Hobson and Klimmek \cite
{hobson-klimmek} (under slightly different conditions). The two
marginal case corresponds to the work of Brown, Hobson and Rogers \cite
{brownhobsonrogers}.

The stochastic control approach allows us to derive the upper bound on
the distribution of $X^*_{t_n}$ in terms of the intermediate
distributions $X_{t_1}, \ldots,\break X_{t_{n-1}}, X_{t_n}$. To show that
the bound is sharp, we need to construct a martingale that fits the
given marginals and attains the bound. To do this we revert to the SEP
methodology.
First we derive the upper bound by taking expectations in a functional
(pathwise) inequality, which in mathematical finance terms has the
interpretation of a semi-static superhedging. This inequality is
``guessed,'' and crucially, it is the stochastic control methodology
that provides candidates for the static and the dynamic components of
the optimal hedge, that is, the different terms in the pathwise inequality.
Once postulated, the inequality is verified independently without any
use of stochastic analysis tools.
Finally, we show the optimality of the upper bound, under some
technical assumptions on the marginals, by establishing that the
solution to the $n$-marginal SEP obtained in Ob\l\'oj and Spoida \cite
{OblSp} achieves equality in our inequality. We note that the idea to
derive martingale inequalities from pathwise inequalities was pivotal
to the pioneering work on robust pricing and hedging of Hobson \cite
{hobson98} and was recently underlined in Acciaio et al. \cite
{acciaio-beiglbock-penkner-schachermayer-temme}.

\subsubsection*{Mathematical finance motivation}
The problem we consider, as described above, has a clear motivation
coming from mathematical finance. The classical framework underpinning
much of the quantitative finance starts by postulating a stochastic
universe $(\Omega, \mathbb{F}, \mathbb{P})$, which is meant to model
a financial environment and capture its riskiness. What it fails to
capture, however, is the uncertainty in the choice of $\mathbb{P}$,
that is, the possibility that the model itself is wrong, also called
\emph{the Knightian uncertainty}; see Knight \cite{knight21}. To
account for model uncertainty it is natural to consider simultaneously
a whole family $\{\mathbb{P}_\alpha\dvtx  \alpha\in\mathcal{A}\}$ of
probability measures. When all $\mathbb{P}_\alpha$ are absolutely
continuous w.r.t. one reference measure $\mathbb{P}$, we speak of
\emph{drift uncertainty} or \emph{dominated} setting. This has
important implications for portfolio choice problems (see F\"ollmer,
Schied and Weber \cite{FSW}) but is not different from an incomplete
market setup in terms of option pricing. However, the \emph
{nondominated} setup when $\mathbb{P}_\alpha$ may be mutually
singular posed new challenges and was investigated starting with
Avellaneda et al. \cite{AvellanedaLevyParas95} and Lyons \cite
{Lyons95}, through Denis and Martini \cite{DM}, to several recent
works, for example, Peng \cite{PengG}, Soner, Touzi and Zhang \cite
{stz-aggregation}, Dolinsky and Soner \cite{DolinskySoner} and
Bouchard and Nutz \cite{BN13}.

Naturally as one relaxes the classical setup, one has to abandon its
precision: under model uncertainty we do not try to have a unique price
but rather to obtain an interval of no-arbitrage prices. Its bounds are
given by seller's and buyer's ``safe'' prices, the superreplication and
the subreplication prices, which can be enforced by trading strategies
that work in all considered models. These bounds can be made more
efficient by enlarging the set of hedging instruments. Indeed, in the
financial markets certain derivatives on the underlying we try to model
are liquid and have well-defined market prices. Without one fixed
model, these options can be included in traded assets without creating
an arbitrage opportunity. By allowing one to trade dynamically in the
underlying and statically (today) in a range of options, one hopes to
have a more efficient approach with smaller intervals of possible
no-arbitrage prices. This constitutes the basis of the so-called \emph
{robust approach to pricing and hedging}.

We contribute to this literature. Our aim was to derive in an explicit
form the superhedging cost of a Lookback option given that the
underlying asset is available for frictionless continuous-time trading,
and that European options for all strikes are available for trading for
a finite set of maturities.
In a zero interest rate financial market, it essentially follows from
the no-arbitrage condition, as observed by Breeden and Litzenberger
\cite{bre}, that these trading possibilities restrict the underlying
asset price process into being a martingale with given marginals. Since
a martingale can be written as a time changed Brownian motion,
and the maximum of a continuous processes is not altered by such a time change,
the one-marginal constraint version of this problem can be converted
into the framework of the Skorokhod embedding problem (SEP). This
observation is the starting point of the seminal paper by Hobson
\cite{hobson98} who exploited the already known optimality result of
the Az\'ema--Yor solution to the SEP and, importantly, provided
an explicit static superhedging strategy. This methodology was
subsequently used to derive robust prices and super/sub-hedging
strategies for barrier options in Brown, Hobson and Rogers \cite
{brownhobsonrogersMF}, for options on local time in Cox, Hobson and
Ob\l\'oj \cite{CoxHobsonObloj08}, for double barrier options in Cox
and Ob\l\'oj \cite{coxobloj,coxobloj2} and for options on variance in
Cox and Wang \cite{coxwang}; see Ob\l\'oj \cite{OblojEQF10} and
Hobson \cite{HobsonSurvey10} for more details.

The above works focused on finding explicitly robust prices and hedges
for an option maturing at $T$ and given market prices of call/put options
co-maturing at~$T$. For lookback options, an extension to the case
where prices at one intermediate maturity are given can be deduced from
Brown, Hobson and Rogers \cite{brownhobsonrogers}. More recently,
Hobson and Neuberger \cite{HobsonNeuberger12} treated forward
starting straddle also using option prices at two maturities.
Otherwise, and excluding the trivial cases when intermediate laws have
no constraining effect (see, e.g., the iterated Az\'ema--Yor setting
in Madan and Yor \cite{mad}), we are not aware of any explicit robust
pricing/hedging results when prices of call options for several
maturities are given. The most likely reason for this is that the
SEP-based methodology pioneered in Hobson \cite{hobson98} starts with
a good guess for the superhedge/embedding, and these become much more
difficult when more marginals are involved. As explained above, our
approach uses stochastic control methods to \emph{derive} a candidate
optimal superhedge strategy. On the dual side, it sees superhedging as
a martingale transportation problem: maximize the expected coupling
defined by the payoff so as to transport the Dirac
measure along the given distributions $\mu_1,\ldots,\mu_n$ by means
of a continuous-time process restricted to be a martingale. This
approach was simultaneously suggested by Beiglb\"{o}ck, Henry-Labord\`
ere and Penkner
\cite{beiglbockhenrylaborderepenkner} in the discrete-time case, and
Galichon, Henry-Labord\`ere and Touzi \cite{ght} in continuous-time.
We refer to Bonnans and Tan \cite{bonnanstan} for a numerical
approximation in the context of variance options, and Tan and Touzi
\cite{TanTouzi} for a general version of the optimal transportation
problem under controlled dynamics.

\subsection*{Organization of the paper}
The paper is organized as
follows. Section~\ref{sect-formulation} provides the precise
mathematical formulation of the problem and establishes the relevant
connections with martingale inequalities, the Skorohod embedding
problem and the martingale optimal transport. The main results of this
paper are collected in Section~\ref{secmainresult}, starting from a
remarkable pathwise inequality. The proofs are reported in Section~\ref
{sectpathwise}. In particular, the pathwise inequality follows from an
elementary verification. The stochastic control approach, which allowed
us to derive the correct quantities for the pathwise arguments, is
pursued in Section~\ref{sectstochasticcontrol}. Additional arguments
for the one marginal case are given in Section~\ref{sectAzema-Yor}.
Finally, the \hyperref[append]{Appendix} contains some proofs of technical lemmas
including some additional properties of the embedding obtained in \cite{OblSp}.

%s2 #&#
\section{Robust superhedging of Lookback options}\label{sect-formulation}

%s2.1 #&#
\subsection{Modeling the volatility uncertainty}\label{sect-setup}

Let $\Omega_x:= \{\omega\in C([0,T], \mathbb{R}^1)\dvtx\break  \omega_0=x\}$.
We consider the set of paths $\Omega:=\Omega_0$ as the canonical
space equipped with the uniform norm $\llVert  \omega\rrVert   _\infty:=
\sup_{0\le t\le T}\llvert  \omega_t\rrvert  $, $B$ the\vspace*{1pt} canonical process, $\mathbb
{P}_0$ the
Wiener measure, $\mathbb{F}:= \{{\mathcal F}_t\}_{0\le t\le T}$ the filtration
generated by $B$. Throughout the paper, $X_0$ is some given initial
value in $\mathbb{R}$, and we denote
\[
X_t := X_0+B_t\qquad\mbox{for } t\in[0,T].
\]
In order to model the volatility uncertainty, we introduce the set
${\mathcal P}$ of all probability measures on $(\Omega,{\mathcal F})$
such that $B$ is a $\mathbb{P}$-martingale. The coordinate process
stands for the price process of an underlying security, and the
restriction to martingale measures ${\mathcal P}$ is motivated by the
classical no-arbitrage results in mathematical finance and will be
justified by the duality results in Theorems~\ref{thm-nmarginal-pricing} and~\ref{thm-nmarginal}.

The quadratic variation process $\langle X\rangle=\langle B\rangle$
is universally defined and takes values in the set of all
nondecreasing continuous functions with $\langle B\rangle_0=0$.

%s2.2 #&#
\subsection{Dynamic trading strategies}\label{sectdynamictrading}

For all $\mathbb{P}\in{\mathcal P}$, we denote by ${\mathbb
{H}}^0(\mathbb{P})$ the collection of all $(\mathbb{P},\mathbb
{F})$-progressively measurable processes and
\begin{eqnarray*}
{\mathbb{H}}^2(\mathbb{P}) &:=& \biggl\{H\in{\mathbb{H}}^0(
\mathbb{P})\dvtx \int_0^T \llvert H_t
\rrvert ^2\,d\langle B\rangle_t<\infty, \mathbb{P}\mbox{-a.s.} \biggr\}.
\end{eqnarray*}
A dynamic trading strategy is defined by a process $H\in\hat{\mathbb
{H}}^2:=\bigcap_{\mathbb{P}\in{\mathcal P}} {\mathbb{H}}^2(\mathbb
{P})$, where $H_t$ denotes the number of shares of the underlying asset
held by the investor at each time $t\in[0,T]$. Under the
self-financing condition, the portfolio value process induced by a
dynamic trading strategy is
%
%e2.1 #&#
\begin{equation}
\label{XH} Y_t^H := Y_0+\int
_0^t H_s \,dB_s,\qquad
t\in[0,T], \mathbb{P}\mbox{-a.s. for all }\mathbb{P}\in{\mathcal P}.
\end{equation}
The stochastic integral in (\ref{XH}) is well defined and should
be rather denoted ${Y^H_t}^\mathbb{P}$ to emphasize its dependence on
$\mathbb{P}$; see, however, Nutz \cite{Nutz}. Nevertheless, for a
large class of strategies $H$ we may define $Y_t^H$ pathwise.
In particular, consider $H\dvtx \Omega_{X_0}\times[0,T]\to\mathbb{R}$ to
be a process of finite variation which is \emph{progressively
measurable} in the sense that $H_t(\omega)=H_t(\omega')$ for any
$t\in[0,T]$ and any $\omega,\omega'\in\Omega_{X_0}$ with $\omega
_s=\omega'_s$, $s\leq t$. We define its
integral (see Dolinsky and Soner \cite{DolinskySoner}) through
an integration by parts formula using classical Stieltjes integration,
%
%e2.2 #&#
\begin{eqnarray}\label{eqpathwiseSI}
\int_0^t H_s(\omega)
\,d\omega_s := H_t(\omega)\omega_t-H_0
\omega _0-\int_0^t
\omega_s \,dH_s(\omega),
\nonumber\\[-8pt]\\[-8pt]
\eqntext{t\in[0,T], \omega\in
\Omega_{X_0}.}
\end{eqnarray}
Note that this integral agrees a.s. with the It\^o stochastic integral
$\mathbb{P}$-a.s. for any $\mathbb{P}\in{\mathcal P}$. We will use
this approach in particular in Section~\ref{secpathwisesuperhedging}
and in Theorem~\ref{thm-nmarginal}.

%s2.3 #&#
\subsection{Semi-static hedging strategies}\label{sectsemistatic}
Let $n$ be some positive integer and $0=t_0<\cdots<t_n=T$ be some
partition of the interval $[0,T]$. In addition to the continuous-time
trading of the primitive securities, we assume that the investor can
take static positions in European call or put options with all possible
strikes and maturities $t_{1} < \cdots<t_{n}$. The market price of the
European call option with strike $K\in\mathbb{R}$ and maturity
$t_{i}$ is denoted
\[
c_{i}(K),\qquad i=1,\ldots,n\quad\mbox{and we denote}\quad c_0(K):=(X_0-K)^+.
\]
A model $\mathbb{P}\in{\mathcal P}$ is said to be calibrated to the
market if $\E^\mathbb{P}[(X_{t_i}-K)^+]=c_i(K)$ for all $1\leq i\leq
n$ and $K\in\mathbb{R}$. For such a model, it was observed by Breeden
and Litzenberger \cite{bre} that, by direct differentiation with
respect to $K$,
\[
\mathbb{P}(X_{t_i}>K)=-c_i'(K+)=:
\mu_{i}\bigl((K,\infty)\bigr),
\]
so that the marginal distributions of $X_{t_i}$, $i=1,\ldots,n$, are
uniquely specified by the market prices and are independent of $\mathbb
{P}$. Let $\mu=(\mu_1,\ldots,\mu_n)$ and
\begin{eqnarray*}
{\mathcal P}^\mu &:=& \{\mathbb{P}\in{\mathcal P}\dvtx  X_{t_i}
\sim\mu_i, 1\le i\le n \}
\end{eqnarray*}
be the set of calibrated market models. By the Strassen theorem \cite
{Strassen}, we have ${\mathcal P}^\mu\neq\varnothing$ if and only if
$\mu_i$'s are nondecreasing in convex order or, equivalently,
%
%e2.3 #&#
\begin{eqnarray}\label{ci-nondecreasing}\label{noarbitrage}
\int\llvert x\rrvert \,d\mu_i(x)<\infty,\qquad
\int x \,d\mu_i(x)=X_0\quad\mbox{and}\quad
c_{{i-1}}\leq c_{i}
\nonumber\\[-8pt]\\[-8pt]
\eqntext{\mbox{for all } 1\leq i\leq n,}
\end{eqnarray}
where now $c_i(K)=\int_{K}^\infty(x-K)\,d\mu_i(x)$. The necessity
follows from Jensen's inequality. For sufficiency, an explicit model
can be constructed using techniques of Skorokhod embeddings; see Ob\l\'
oj \cite{Obloj04b}. Consequently, up to integrability, the
$t_i$-maturity European derivative defined by the payoff $\lambda
_{i}(X_{t_{i}})$ has an unambiguous market price
\begin{eqnarray*}
&&\mu_{i}(\lambda_{i}) := \int\lambda_{i} \,d
\mu_{i} =\E^\mathbb{P}\bigl[\lambda_i(X_{t_i})
\bigr]\qquad\mbox{for all } \mathbb{P}\in{\mathcal P}^\mu.
\end{eqnarray*}
The condition ${\mathcal P}^\mu\neq\varnothing$ embodies the fact that
the current observed market prices do not induce arbitrage. By this we
mean that there exists a model which admits no-arbitrage (no free lunch
with vanishing risk) and reprices the call options through risk neutral
expectation. For this reason we sometimes refer to {(\ref
{noarbitrage})} as the \emph{no-arbitrage condition}. We note,
however, that the arbitrage considerations are more subtle in all
generality since boundary cases may arise where ${\mathcal P}^\mu=
\varnothing$ but there is no model-independent arbitrage; see Davis and
Hobson \cite{davishobson} and Cox and Ob\l\'oj \cite{coxobloj}.

%
%re2.1 #&#
\begin{Remark}
For the purpose of the present financial application, we could
restrict the measures $\mu_i$ to have support in $\mathbb{R}_+$ and
$\mathbb{P}\in{\mathcal P}$ to be such that $X_t\geq0$ $\mathbb
{P}$-a.s. Note, however, that this is easily achieved: it suffices to
assume that \mbox{$X_0>0$} and $c_n(K)=X_0-K$ for $K\leq0$. Then $\mu
_n((K,\infty))=1$, $K<0$, and hence $\mu_n([0,\infty))=1$. Then for
any $\mathbb{P}\in{\mathcal P}^\mu$ we have $X_t=\E^{\mathbb
{P}}[X_T\mid \mathcal{F}_t]\geq0$ $\mathbb{P}$-a.s. for $t\in[0,T]$.
In particular, $\mu_i([0,\infty))=\mathbb{P}(X_{t_i}\geq0)=1$.
\end{Remark}

We denote $\tbf:=(t_{1},\ldots,t_{n})$,
$\lambda=(\lambda_{1},\ldots,\lambda_{n})$,
%
%e2.4 #&#
\begin{equation}
\label{Glambda} \mu(\lambda):=\sum_{i=1}^n
\mu_{i}(\lambda_{i}), \qquad \lambda(\omega_\tbf) :=
\sum_{i=1}^n\lambda_{i}(
\omega_{t_{i}}),
\end{equation}
for $\omega\in C([0,T])$.
The set of Vanilla payoffs which may be used by the hedger are
naturally taken in the set
%
%e2.5 #&#
\begin{eqnarray}
\label{Lambdan} \Lambda^\mu_n &:=& \biggl\{\lambda\dvtx \int
\llvert \lambda_{i}\rrvert \,d\mu_i<\infty, 1\le i\le n
\biggr\}.
\end{eqnarray}
Thus, in addition to the dynamic hedging strategies $H$, the investor
has access to the static hedging strategies $\lambda$, consisting of a
buy-and-hold strategy in a portfolio of options. Such a pair $(\lambda
,H)$ is called a semi-static hedging strategy and induces the final
value of the self-financing portfolio
%
%e2.6 #&#
\begin{eqnarray}
\label{barX} \overline{Y}{}^{H,\lambda}_T &:=&
Y^H_T-\mu(\lambda)+\lambda(X_\tbf),
\end{eqnarray}
indicating that the investor has the possibility of buying at time $0$
any derivative security with payoff $\lambda_i(X_{t_i})$ for the price
$\mu_i(\lambda_i)$. As it will be made clear in our subsequent
analysis, the functions $\lambda_{i}$ will play the role of a Lagrange
multiplier for the constraints $X_{t_{i}}\sim\mu_{i}$,
$i=1,\ldots,n$.

%s2.4 #&#
\subsection{Robust superhedging under semi-static hedging strategies}

In this paper, we focus on the problem of robust semi-static
superhedging of a lookback option defined by the payoff at the maturity $T$,
\begin{eqnarray*}
&&\textstyle \xi:= \phi (X^*_T)\mbox{ where } X^*_t:=\max
_{s\le t}X_s\mbox{ and }\phi\dvtx \mathbb{R}\longmapsto
\mathbb{R}\mbox{ is right-continuous,}
\\[-4pt]
&&\mbox{nondecreasing.}
\end{eqnarray*}
The investor can trade as discussed in the previous two sections.
However, we need to impose a further admissibility condition to rule
out doubling strategies. We let ${\mathcal H}^\mu$ consist of all
processes $H\in\hat{\mathbb{H}}^2$ whose induced portfolio value
process $Y^H$ is $\mathbb{P}$-supermartingale for all $\mathbb{P}\in
{\mathcal P}^\mu$.
The robust superhedging upper bound is then defined by
%
%e2.7 #&#
\begin{equation}
\label{V}
\qquad U^\mu_n(\xi) := \inf \bigl\{
Y_0\dvtx  \exists (\lambda,H)\in\Lambda^\mu_n\times
{\mathcal H}^\mu, \overline{Y}{}^{H,\lambda}_T\ge\xi,
\mathbb{P}\mbox{-a.s. for all }\mathbb{P}\in{\mathcal P} \bigr\}.
\end{equation}
Selling $\xi$ at a price higher than $U^\mu_n(\xi)$, the hedger
could set up a portfolio with a negative initial cost and a
nonnegative payoff under any market scenario leading to a strong
(model-independent) arbitrage opportunity. Note that, thanks to the
definition of admissible trading strategies, by taking expectations in
the superhedging inequality $\overline{Y}{}^{H,\lambda}_T\ge\xi$, and
optimizing over $(\lambda,H)$ and $\mathbb{P}$, we obtain the usual
pricing-hedging inequality
%
%e2.8 #&#
\begin{equation}
\label{eqintrophineq} U^\mu_n(\xi)\geq\sup_{\mathbb{P}\in{\mathcal P}^\mu}
\E^\mathbb {P}[\xi].
\end{equation}
Theorems~\ref{thm-nmarginal-pricing} and~\ref{thm-nmarginal} below
establish a bound on the RHS value and show that under mild technical
conditions equality holds in {(\ref{eqintrophineq})} and exhibit
both the best superhedge and the maximal $\mathbb{P}^{\max}$
on the RHS. However, before proceeding to our main results, we first
discuss the connection of the superhedging problem with two important
questions in applied mathematics, namely martingale inequalities in
probability and the theory of optimal transport.

%s2.5 #&#
\subsection{Pathwise super-hedging and martingale inequalities}\label{secpathwisesuperhedging}

In this section, we discuss briefly the connection between the robust
super-hedging problem and pathwise inequalities inducing martingale
inequalities, as highlighted by Acciaio et~al.~\cite{acciaio-beiglbock-penkner-schachermayer-temme}
and Beiglb\"ock and Nutz \cite{BN14}; see also
Os\c{e}kowski \cite{O12}.
Suppose that existence holds in the super hedging problem {(\ref
{V})}. Then $\overline{Y}{}^{H,\lambda}_0=U^\mu_n(\xi)$ and
$\overline{Y}{}^{H,\lambda}_T\ge\xi$, $\mathbb{P}$-a.s. for all
$\mathbb{P}\in{\mathcal P}$, for some $\lambda\in\Lambda^\mu_n$
and $H\in{\mathcal H}^\mu$. Assume further that the dynamic hedging
process $H$ is such that its ``stochastic integral'' can be defined
pathwise and that the super-hedging property extends to a pathwise inequality
%
%e2.9 #&#
\begin{equation}
\label{eqpathwisesuperhedge}
\qquad \xi(\omega) \le U^\mu_n(\xi) +\lambda(
\omega_{\mathbf{t}})-\mu(\lambda) +\int_0^T
H_s \,d\omega(s)\qquad \mbox{for all } \omega\in\Omega_{X_0}.
\end{equation}
This inequality is then sharp in the sense that the above pair
$(\lambda,H)$ minimizes the cost $\mu(\lambda)$ of the trading
strategy among all such strategies which super-hedge $\xi$ pathwise.
In particular, $U^\mu_n(\xi)$ is the superhedging price not only in
the sense of~{(\ref{V})} but also in the pathwise sense.
Assuming further that the stochastic integral in~{(\ref
{eqpathwisesuperhedge})} defines a uniformly integrable martingale,
this implies a martingale inequality
\begin{eqnarray*}
&& \E[\xi(Y)] \le U^{{\mathcal L}^Y_{\mathbf{t}}}_n(\xi)
\mbox{ for any continuous martingale } (Y_t)_{t\leq T},
\\[-4pt]
&&\mbox{where }{\mathcal
L}^Y_{\mathbf{t}}=({\mathcal L}^{Y_{t_i}}, i\le n),
\end{eqnarray*}
and ${\mathcal L}^{Y_{t_i}}$ is the distribution of $Y_{t_i}$ and
$Y_0=X_0$. Moreover, by construction, this inequality is sharp: we can
construct continuous martingales $Y$ which attain equality.
An example of such martingale inequality is provided in Proposition
\ref{propmartineq} below.

%s2.6 #&#
\subsection{Optimal transportation under martingale restriction}\label{secrobust-SEP}

In this short section we discuss the connection of our problem to
optimal transportation theory, which will be the building block for the
stochastic control approach of Section~\ref{sectstochasticcontrol}.

The first duality we consider is the expected quasi-sure extension of the
classical dual formulation of the superhedging problem
\begin{eqnarray*}
U^\mu_n(\xi) &=& \inf_{\lambda\in\Lambda^\mu_n} \sup
_{\mathbb{P}\in{\mathcal P}} \E^\mathbb{P} \bigl[\xi+\lambda(X_{\mathbf{t}})-
\mu(\lambda) \bigr].
\end{eqnarray*}
We show it holds, subject to
technical assumptions, in Proposition \ref{propduality} below.
The second duality follows by formally inverting the inf-sup on the RHS
above leading to the optimization problem
%
%e2.10 #&#
\begin{equation}
\label{transport} \sup_{\mathbb{P}\in{\mathcal P}^\mu} \E^\mathbb{P}[\xi],
\end{equation}
which falls into the recently introduced class of optimal
transportation problems under controlled stochastic dynamics; see
Beiglb\"{o}ck, Henry-Labord\`ere and Penkner \cite
{beiglbockhenrylaborderepenkner}, Galichon, Henry-Labord\`ere and Touzi
\cite{ght} and Tan and Touzi \cite{TanTouzi}. In words, the
above problem consists of maximizing the expected transportation
cost of the Dirac measure $\delta_{\{X_0\}}$ along the given
marginals $\mu_1,\ldots,\mu_n$ with transportation scheme
constrained to the class of martingales. The cost of transportation in
our context is defined by the path-dependent payoff $\xi(\omega)$.

The validity of the equality between the value function in
{(\ref{transport})} and our problem $U^\mu_n(\xi)$ was
established recently by Dolinsky and Soner \cite{DolinskySoner} for
$n=1$ and under strong continuity assumptions on the payoff function
$\omega\longmapsto\xi(\omega)$. The corresponding duality result in
the discrete time framework was obtained by Beiglb\"{o}ck,
Henry-Labord\`ere and Penkner \cite{beiglbockhenrylaborderepenkner}.

Note that if we can find a trading strategy $\overline{Y}{}^{H,\lambda
}_T$ as in {(\ref{barX})} which superreplicates $\xi$: $\overline
{Y}{}^{H,\lambda}_T\geq\xi$ $\mathbb{P}$-a.s. for all $\mathbb{P}\in
{\mathcal P}$ and a $\mathbb{P}^{\max}\in{\mathcal P}^\mu$ such
that $\E^{\mathbb{P}^{\max}}[\xi]=Y_0$, then as in {(\ref
{eqintrophineq})},
\[
Y_0 \le \sup_{\mathbb{P}\in{\mathcal P}^\mu}\E^\mathbb{P}[\xi] \le
U^\mu_n(\xi) \le Y_0,
\]
and it follows that we have equalities throughout. This line of attack
has been at the heart of the approach to robust pricing and hedging
based on the Skorokhod embedding problem, as in
Hobson \cite{hobson98}, Brown, Hobson and Rogers \cite
{brownhobsonrogersMF}, Cox and Ob{\l}{\'o}j \cite{coxobloj,coxobloj2}
and Cox and Wang \cite{coxwang}.
It relies crucially on the ability to make a correct guess for the
cheapest superhedge $\overline{Y}{}^{H,\lambda}_T$. This becomes
increasingly difficult when one considers information about prices at
several maturities, $n>1$. In this paper, we follow the above
methodology in Section~\ref{sectpathwise} to prove our main result,
Theorem~\ref{thm-nmarginal}.
Sections~\ref{sectstochasticcontrol}--\ref{sectAzema-Yor} then
provide an alternative approach based on stochastic control methods.
The latter is longer and more involved than the former and requires
slight modifications for technical reasons. However, it is in fact
necessary as it is the origin of the determination of the right
quantities for the former, namely the pathwise inequality.

%s2.7 #&#
\subsection{The Skorokhod embedding problem}\label{secSEP}

We now specialize the discussion to the case of a lookback option $\xi
=G(X_\tbf,X^*_T)$, for some payoff function $G$. By the
Dambis--Dubins--Schwarz theorem, we may
re-write problem {(\ref{transport})} as a multiple stopping problem
(see Proposition 3.1 of Galichon, Henry-Labord\`ere and Touzi~\cite{ght}),
%
%e2.11 #&#
\begin{equation}
\label{SEP} \sup_{(\tau_1,\ldots,\tau_n)\in{\mathcal T}^\mu} \E^{\mathbb{P}_0} \bigl[G
\bigl(X_{\tau_1},\ldots,X_{\tau
_n},X^*_{\tau_n} \bigr)
\bigr],
\end{equation}
where the ${\mathcal T}^\mu$ is the set of ordered stopping times
$\tau_1\leq\cdots\leq\tau_n<\infty$ $\mathbb{P}_0$-a.s. with $X_{\tau_i}\sim_{\mathbb{P}_0}\mu_i$ for all $i=1,\ldots,n$
and $(X_{t\land\tau_n})$ being a uniformly integrable martingale.
Elements of ${\mathcal T}^\mu$ are solutions to the iterated
(multi-marginal) version of the so-called Skorokhod
embedding problem (SEP); cf. Ob{\l}{\'o}j \cite{Obloj04b}.
Here, formulation {(\ref{SEP})} is directly searching for a
solution to the
SEP which maximizes the criterion defined by the coupling $G(x,m)$.
Previous works have focused mainly on single marginal constraint ($n=1$).
The case $G(x,m)=\phi(m)$ for some nondecreasing function $\phi$ is
solved by the so-called Az\'ema--Yor embedding; cf. Az\'ema and Yor
\cite{AzemaYor1,AzemaYor2}, Hobson \cite{hobson98}; see also
Galichon, Henry-Labord\`ere and Touzi \cite{ght} who recovered this
result by the stochastic-control approach of Section~\ref
{sectstochasticcontrol}. The case $G(x,m)$ was considered
recently by Hobson and Klimmek \cite{hobson-klimmek}, where the
optimality of the Az\'ema--Yor solution of the SEP is shown to be
valid under convenient conditions on the function $G$. This case is
also solved in Section~\ref{sectAzema-Yor} of the present paper, with
our approach leading to the same results as those obtained by Hobson
and Klimmek \cite{hobson-klimmek}, but under slightly different assumptions.

The case $G(x_1,\ldots,x_n,m)=\phi(m)$ for some nondecreasing
function $\phi$ is also trivially solved by $\tau^{AY}(\mu_n)$ in
the following special case when the single marginal solutions are
naturally ordered: $\tau^{AY}(\mu_i)\leq\tau^{AY}(\mu_{i+1})$.
This is called the \textit{increasing mean residual value
property} by Madan and Yor \cite{mad} who establish, in particular, a
strong Markov property of the resulting time-changed process.
The case of arbitrary measures which satisfy {(\ref{noarbitrage})}
for $n=2$ was solved in Brown, Hobson and Rogers \cite{brownhobsonrogers}.
In this paper we consider $n\in\mathbb{N}$. For subsequent use, we
recall that the Az\'ema--Yor embedding for $\mu_i$ is given by $\tau
^{AY}(\mu_i)=\inf\{t\geq0\dvtx  X_t\leq b_i^{-1}(X_t^*)\}$, where the
inverse barycentre function $b^{-1}_i(m)$ is a minimizer in
%
%e2.12 #&#
\begin{equation}
\label{eqone-dimoptimisation} %\frac{c_i\left(b^{-1}_i(m)\right)}{m-b^{-1}_i(m)}=
\min_{\zeta\leq m} \frac{c_i(\zeta)}{m-\zeta},\qquad
m>X_0,
\end{equation}
taken to be right-continuous in $m$ and with $b^{-1}_i(m)=m$ for $m\geq
\inf\{m\dvtx\break  c_i(m)=0\}$. It is easy to see that $b^{-1}_i$ is
nondecreasing. Note also that $c_i(\zeta)/(m-\zeta)$ is
nonincreasing for $\zeta\leq b^{-1}_i(m)$ and nondecreasing for
$\zeta\geq b^{-1}_i(m)$.

%s3 #&#
\section{Main results}\label{secmainresult}

Our main result is split into three parts. We first state a
trajectorial inequality which is the building block for the solution of
the robust superhedging problem. We next solve the pricing problem {(\ref{transport})} in Theorem~\ref{thm-nmarginal-pricing}. Finally,
in Theorem~\ref{thm-nmarginal}, we solve the superhedging problem {(\ref{V})}.

%s3.1 #&#
\subsection{A remarkable trajectorial inequality}

The first result involves, for all $\zeta_1 \leq\cdots\leq\zeta_n <
m$, the semi-static hedging strategy
%
%e3.1 #&#
%e3.2 #&#
\begin{eqnarray}
\lambda_i^{\zeta,m}(x) & := & \frac{(x-\zeta_i)^+}{m-\zeta_i}-{
\mathbf1}_{ \{ {i<n}  \}
}\frac{(x-\zeta_{i+1})^+}{m-\zeta_{i+1}}, \qquad x\in\mathbb{R},
\\
H^{\zeta,m}_t(\omega) & := & -\frac{{\mathbf1}_{(t_{i-1},t]}(T_m(\omega))
+{\mathbf1}_{[0,t_{i-1}]}(T_m(\omega)){\mathbf1}_{ \{ {\omega
_{t_{i-1}}\ge\zeta_i}  \} }}{
m-\zeta_i},
\nonumber\\[-8pt]\\[-8pt]
\eqntext{t \in[t_{i-1},t_i),}
\end{eqnarray}
for all $i=1,\ldots,n$, where $T_m(\omega):=\inf\{t\ge0\dvtx \omega
_t\ge m\}$. Notice that $H^{\zeta,m}$ is a piecewise constant
predictable process, so that the stochastic integral $\int H_t(\omega
)\,d\omega_t$ is a well-defined, pathwise, finite sum.

%
%pr3.1 #&#
\begin{Proposition}
\label{lemTrajectorialInequalityOrderedCase}
Let $\omega$ be a c\`{a}dl\`{a}g path, $m>\omega_0$, and denote
$\omega^*_t:=\sup_{0\le s \le t}{\omega_s}$. Then, for all $\zeta_1
\leq\cdots\leq\zeta_n < m$, the following inequality holds:
%
%e3.3 #&#
\begin{eqnarray}
{\mathbf1}_{ \{ {\omega^*_{t_n} \geq m}  \} } &\le& \sum_{i=1}^{n}
\biggl\{\lambda_i^{\zeta,m}(\omega_{t_i}) +\int
_{t_{i-1}}^{t_i} H^{\zeta,m}_t(\omega) \,d
\omega_t \biggr\}. \label{eqthmTrajectorialInequalityOrderedCase1}
\end{eqnarray}
\end{Proposition}

The proof is reported in Section~\ref{secttrajectorial} and is
based on an elementary verification. The importance of this inequality
is that it will be shown to be sharp in some precise sense, so that the
solution of the robust superhedging problem is fully deduced from it.
Thereore, the relevant difficulty is on how this inequality can be
guessed. This issue is addressed in Section~\ref
{sectstochasticcontrol}, where our intention is to show that the
stochastic control approach is genuinely designed for this purpose.

The pathwise inequality of Proposition~\ref
{lemTrajectorialInequalityOrderedCase} is stated for the elementary
lookback option defined by the payoff function $\phi={\mathbf
1}_{[m,\infty)}$. The corresponding extension to a general
right-continuous, nondecreasing function $\phi$ follows by the
obvious identity
%
%e3.4 #&#
\begin{eqnarray}
\label{indicatortophi} \phi\bigl(\omega^*_{t_n}\bigr) &=& \phi(
\omega_0)+\int_{(\omega_0,\infty)} {\mathbf1}_{ \{ {\omega
^*_{t_n} \geq m}  \} }\,d
\phi(m).
\end{eqnarray}

%s3.2 #&#
\subsection{Financial interpretation}\label{subsecFinancialInterpretation}

We develop now a financial interpretation of the RHS of
{(\ref{eqthmTrajectorialInequalityOrderedCase1})}, for
$\omega=X$ the price process, as a (pathwise) superhedging strategy
for a simple knock-in digital barrier option with payoff\vspace*{2pt} $\xi={\mathbf
1}_{ \{ {X^*_T \geq m}  \} }$. The semi-static hedging
strategy consists of three elements: a static position in call options,
a forward transaction (with the shortest available maturity) when the
barrier $m$ is hit and rebalancing thereafter at times $t_i$.
More precisely:

\begin{longlist}[(iii)]
\item[(i)] \textit{Static position in calls}:
\begin{eqnarray*}
\lambda^{\zeta,m}(X_\tbf) &:=& \sum
_{i=1}^{n} \lambda^{\zeta,m}_i(X_{t_i}).
\end{eqnarray*}
For $1\leq i<n$, we hold a portfolio long and short calls with maturity
$t_i$ and strikes $\zeta_i$ and $\zeta_{i+1}$, respectively. This
yields a ``tent like'' payoff which becomes negative only if the
underlying exceeds level $m$. Note that by setting $\zeta_i=\zeta
_{i+1}$ we may avoid trading the $t_i$-maturity calls. For maturity
$t_n$ we are only long in a call with strike $\zeta_n$.

\item[(ii)] \textit{Forward transaction if the barrier $m$ is hit}:
\[
\int_{t_{i-1}}^{t_{i}} H^{\zeta,m}_t(X)
\,dX_t = \frac
{m-X_{t_i}}{m-\zeta_i}\qquad\mbox{on }\bigl\{t_{i-1}<T_m(X)
\leq t_i\bigr\}=\bigl\{ X^*_{t_{i-1}}<m\leq X^*_{t_i}
\bigr\}.
\]
At the moment when the barrier $m$ is hit,\footnote{This is well
defined since we only consider continuous paths. Note, however, that
even if the process was allowed to jump over the level $m$, the
superhedging property would be preserved and only an additional profit
would be realized.}
say between maturities $t_{i-1}$ and $t_{i}$, we enter into forward
contracts with maturity $t_{i}$.
Note that the long call position with maturity $t_{i}$ together with
the forward then superhedge the knock-in digital barrier option; cf.
{(\ref{eqforwardtransaction})}. This resembles the robust
semi-static hedge in the one-marginal case; cf. Lemma 2.4 of Brown,
Hobson and Rogers \cite{brownhobsonrogers}. All the ``tent like''
payoffs up to maturity $t_{i-1}$ are nonnegative.

\item[(iii)] \textit{Rebalancing of portfolio to hedge calendar spreads}:
\[
\int_{t_i}^{t_n} H^{\zeta,m}_t
\,dX_t = \sum_{j=i}^n {
\mathbf1}_{ \{
{X_{t_j}\geq\zeta_{j+1}}  \} }\frac
{X_{t_j}-X_{t_{j+1}}}{m-\zeta_{j+1}}\qquad\mbox{on }\bigl\{ t_{i-1}<T_m(X)
\leq t_i\bigr\}.
\]
After the barrier $m$ is hit between $t_{i-1}$ and $t_i$, we start
trading at times $t_{j}$, $j\geq i$, in such a way that a potential
negative payoff of the calendar spreads $\frac{(X_{t_{j}}-\zeta
_{j})^+}{m-\zeta_{j}} - \frac{(X_{t_{j-1}}-\zeta_{j})^+}{m-\zeta
_{j}}, i < j \leq n$, is offset; cf. {(\ref
{eqrebalancingcalendarspreads})}.

In the above, (ii) and (iii) are instances of dynamic trading which is
done in a self-financing way. Their combined payoff is $\int_0^T
H^{\zeta,m}_s d X_s$, and inequality {(\ref
{eqthmTrajectorialInequalityOrderedCase1})} simply says that for
any choice of $\zeta_1\leq\cdots\leq\zeta_n<m$, the semi-static
hedging strategy $(H^{\zeta,m},\lambda^{\zeta,m})$ superreplicates
$\xi$.
\end{longlist}

%s3.3 #&#
\subsection{Martingale inequalities and robust superhedging}

As a first consequence of the trajectorial inequality of Proposition
\ref{lemTrajectorialInequalityOrderedCase}, we have the following
martingale inequality involving finitely-many intermediate marginals.

%
%pr3.2 #&#
\begin{Proposition}
\label{propmartineq}
Let $Y$ be a c\`adl\`ag submartingale defined on a filtered probability
space satisfying the usual conditions, and consider an arbitrary
right-continuous nondecreasing function $\phi$. Then, for any
functions $\zeta_1(m)\le\cdots\le\zeta_n(m) \leq m$, we have
\begin{eqnarray*}
\E \bigl[\phi\bigl(Y^*_{t_n}\bigr) \bigr] &\le& \phi(Y_0) +
\int_{(X_0,\infty)} \sum_{i=1}^n
\biggl(\frac{\E[(Y_{t_i}-\zeta_i(m))^+]}{
m-\zeta_i(m)}
\\
&&\hspace*{93pt}{} -\frac{\E[(Y_{t_i}-\zeta_{i+1}(m))^+]}{
m-\zeta_{i+1}(m)}{\mathbf1}_{\{i<n\}} \biggr) \,d
\phi(m),
\end{eqnarray*}
where $Y^*_t:= \sup_{u\leq t} Y_u$.
\end{Proposition}

\begin{pf}
Taking expectation in inequality {(\ref
{eqthmTrajectorialInequalityOrderedCase1})} for any $\zeta_1\le
\cdots\leq\zeta_n < m$, we see that
\begin{eqnarray*}
\E [{\mathbf1}_{\{Y^*_{t_n}\ge m\}} ]
&\le&  \sum_{i=1}^n \biggl\{ \E \bigl[
\lambda^{\zeta,m}(Y_{t_i}) \bigr]
\\
&&\hspace*{17pt}{}-\frac{\E [Y_{t_i\vee T_m}-Y_{t_{i-1}\vee T_m}
+{\mathbf1}_{ \{ {T_m \leq t_{i-1},  Y_{t_{i-1}} \geq\zeta_i}
 \} } (Y_{t_i}-Y_{t_{i-1}} )
 ]}{
m-\zeta_i} \biggr\}
\\
&\le&  \sum_{i=1}^n \E \bigl[
\lambda^{\zeta,m}(Y_{t_i}) \bigr],
\end{eqnarray*}
by the submartingale property of $Y$. Taking limits, the above extends
to any $\zeta_1\le\cdots\leq\zeta_n \leq m$ giving the
required inequality for $\phi={\mathbf1}_{[m,\infty)}$. Then, for a
right-continuous nondecreasing function $\phi$, we take expectation
in {(\ref{indicatortophi})}, and we conclude using Fubini and
the last inequality.
\end{pf}

The particular case $\phi={\mathbf1}_{[m,\infty)}$ provides an upper
bound on $\mathbb{P}[Y^*_{t_n}\geq m]$. Note also that for some $Y$,
the RHS may reduce to a much simpler form. In particular if $Y$ is
stopped at $t_1$, $Y_t=Y_{t\land t_1}$ for $t\geq0$, then the RHS
reduces to the one marginal case. We explore martingale inequalities of
the above form and their usefulness in a short parallel note \cite{OST13}.

The key ingredient for the solution of the present $n$-marginals robust
superhedging problem is to re-write the upper bound in the last
martingale inequality in terms of the call prices
\begin{eqnarray*}
\sum_{i=1}^n \E^\mathbb{P} \bigl[
\lambda^{\zeta,m}(X_{t_i}) \bigr] &=& \sum
_{i=1}^n \biggl(\frac{c_i(\zeta_i)}{
m-\zeta_i} -
\frac{c_i(\zeta_{i+1})}{
m-\zeta_{i+1}}{\mathbf1}_{\{i<n\}} \biggr)\qquad\mbox{for all } \mathbb{P}\in{
\mathcal P}^\mu.
\end{eqnarray*}
By the arbitrariness of the parameters $\zeta$, we are then reduced to
the best upper bound
%
%e3.5 #&#
\begin{eqnarray}\label{eqglobaloptimizationproblem}
C(m) &:=& \min_{\zeta_1\le\cdots\le\zeta_n\le m} \sum
_{i=1}^n \biggl(\frac{c_i(\zeta_i)}{
m-\zeta_i} -
\frac{c_i(\zeta_{i+1})}{
m-\zeta_{i+1}}{\mathbf1}_{\{i<n\}} \biggr)
\nonumber\\[-10pt]\\[-8pt]
\eqntext{\mbox{for all } m >
X_0,}
\end{eqnarray}
where here, and throughout, we understand the value of the sum on the
RHS of~{(\ref{eqglobaloptimizationproblem})} for $\zeta
_k<\zeta_{k+1}=\cdots=\zeta_n=m$ as limit of the value $\zeta
_{k+1}=\cdots=\zeta_n=\zeta\to m$ which is clearly either $+\infty$
or is well defined in terms of the derivative of the call function at $m$.

%
%th3.3 #&#
\begin{Theorem}\label{thm-nmarginal-pricing}
Let $\phi$ be a right-continuous nondecreasing function.

\begin{longlist}[(ii)]
\item[(i)] Under the no-arbitrage condition {(\ref
{ci-nondecreasing})}, we have
%
%e3.6 #&#
\begin{eqnarray}
\label{eqlookbackotbound} \sup_{\mathbb{P}\in{\mathcal P}^\mu} \E^\mathbb{P}\bigl[\phi
\bigl(X^*_T\bigr)\bigr] &\le& \phi(X_0)+ \int
_{(X_0,\infty)} C(m)\,d\phi(m).
\end{eqnarray}

\item[(ii)] If in addition $\mu_1,\ldots,\mu_n$ satisfy Assumption
$\circledast$ of Ob\l\'oj and Spoida \cite{OblSp}, then
equality holds in {(\ref{eqlookbackotbound})} and is attained
by some $\mathbb{P}^{\max}\in{\mathcal P}^\mu$.
\end{longlist}
\end{Theorem}

Part (i) of the last theorem is a direct consequence of Proposition
\ref{propmartineq}. The proof of part~(ii) is reported in
Section~\ref{secproofdirect}.
The upper bound {(\ref{eqlookbackotbound})} holds in all
generality; in particular, both sides could be infinite.

We next focus on the existence of a semi-static superhedging strategy
which induces upper bound {(\ref{eqlookbackotbound})}. For the
following preliminary result, we recall the inverse barycenter
functions $b_i$ introduced by {(\ref{eqone-dimoptimisation})}.

%
%le3.4 #&#
\begin{Lemma}\label{27565734534563}
There exists a measurable minimiser $\zeta^*(m)$ for the optimization
problem {(\ref{eqglobaloptimizationproblem})} with $\zeta
^*_1(m)\geq\ub^{-1}(m):= \min_{1\leq i\leq n}b^{-1}_i(m)$.
\end{Lemma}

The proof of this lemma is reported in Appendix~\ref{sectAppendixAA}.
Given these minimizing functions $m\longmapsto\zeta^*(m)$, we deduce
from Proposition~\ref{lemTrajectorialInequalityOrderedCase}
together with {(\ref{indicatortophi})} the following candidates
for the optimal semi-static hedging strategies:
%
%e3.7 #&#
\begin{eqnarray}\label{optimalsemistatic}
\hat\lambda_i &:=& \int_{(X_0,\infty)}
\lambda_i^{\zeta^*_i(m),m}\,d\phi (m),\qquad i=1,\ldots,n,
\nonumber\\[-8pt]\\[-8pt]\nonumber
\hat H_t &:=& \int_{(X_0,\infty)} H_t^{\zeta^*(m),m}\,d
\phi(m),
\end{eqnarray}
where we will impose further assumptions on $\phi$ under which these
integrals are well defined.
Note that from the definition of $T_m$ it follows that
%
%e3.8 #&#
\begin{eqnarray}\label{eqHdecomposition}
-\hat H_t &=& \int_{(\omega^*_{t_{i-1}},\omega^*_t]}
\frac{d\phi(m)}{m-\zeta^*_i(m)}
\nonumber\\[-8pt]\\[-8pt]\nonumber
&&{} + \int_{(X_0,\omega^*_{t_{i-1}}]} {\mathbf1}_{ \{ {\omega_{t_{i-1}}\ge\zeta^*_i(m)}  \} }
\frac{d\phi(m)}{m-\zeta^*_i(m)}, \qquad t\in[t_{i-1},t_i).
\end{eqnarray}
The first term is nondecreasing in $t\in[t_{i-1},t_i)$ while the
second one is constant. It follows that $\hat H$ is of finite
variation. It is also clear that $\hat H$ is progressively measurable
in the sense introduced in Section~\ref{sectdynamictrading} and
hence the stochastic integral $\int\hat H_t \,d\omega_t$ is well-defined pathwise for $\omega\in\Omega_{X_0}$ via {(\ref
{eqpathwiseSI})}.

% Further, we will show below that under
%\begin{equation}\label{eqnodoublingsufficient}
%\int_{\R}\int_{(X_0,\infty)}\indicator{x\geq\zeta^*_{i+1}(m)}\frac{d
%\phi(m)}{m-\zeta_{i+1}^*(m)}\,d\phi(m)\mu_i(\,dx)<\infty,\quad
%i=1,\ldots, n-1,
%\end{equation}
%$\hat H\in\Hc^\mu$ and hence $(\hat\lambda, \hat H)$ is an admissible
%semi-static hedging strategy.
We now show that, under mild technical assumptions, equality holds in
{(\ref{eqintrophineq})} and that $(\hat\lambda, \hat H)$ is
optimal and attains the minimal superhedging cost.

%
%th3.5 #&#
\begin{Theorem}\label{thm-nmarginal}
Assume that the no-arbitrage condition {(\ref{ci-nondecreasing})}
holds, and let $\zeta^{*}_1(m),\ldots,\zeta^{*}_n(m)$ be defined by
Lemma~\ref{27565734534563}. Let $\phi$ be a right-continuous nondecreasing
function such that $m \mapsto (m - \zeta^*_n(m) )^{-1}$ is
$d\phi$-locally integrable, $\hat\lambda$ and $\hat H$ be given by
{(\ref{optimalsemistatic})}, and assume that
%
%e3.9 #&#
\begin{equation}
\label{eqfinitelambdacondition} \int_{\mathbb{R}}\int_{(X_0,\infty)}
\frac{(x-\zeta
_i^*(m))^+}{m-\zeta_i^*(m)}\,d\phi(m)\mu_i(dx)<\infty,\qquad i=1,\ldots, n,
\end{equation}
which is equivalent to $\hat\lambda\in\Lambda_n^\mu$. Then $\hat
H\in{\mathcal H}^\mu$ and:

\begin{longlist}[(ii)]
\item[(i)] $U^\mu_n (\phi(X^*_T) )
 \le
\phi(X_0)+ \mu(\hat{\lambda})
 =  \phi(X_0)+
\int_{(X_0,\infty)}C(m) \,d \phi(m)<\infty$, and
%
%e3.10 #&#
\begin{equation}
\label{hedgingstrategy} \phi\bigl(\omega^*_T\bigr) \leq \phi(X_0)+
\hat\lambda(\omega_{\tbf}) +\int_{0}^T
\hat{H}_t(\omega) \,d \omega_t\qquad\mbox{for all } \omega\in
\Omega_{X_0};
\end{equation}

\item[(ii)] if, in addition, $(\mu_i)_{1\le i\le n}$ satisfy Assumption
$\circledast$ of Ob\l\'oj and Spoida \cite{OblSp}, then
$U^\mu_n (\phi(X^*_T) ) = \phi(X_0)+ \mu(\hat{\lambda
})$, and equality holds in {(\ref{hedgingstrategy})}
$\mathbb{P}^{\max}$-a.s., for $\mathbb{P}^{\max}\in{\mathcal
P}^\mu$ described in the proof.
\end{longlist}
\end{Theorem}

The proof is reported in Section~\ref{secproofdirect1}. Note that
when $\hat\lambda\in\Lambda_n^\mu$ and $\hat H\in{\mathcal H}^\mu
$, claim~(i)~is a direct consequence of the pathwise inequality of
Proposition~\ref{lemTrajectorialInequalityOrderedCase} together
with {(\ref{indicatortophi})} and the assumed conditions.

%
%re3.6 #&#
\begin{Remark}
It follows from Section~4 of Ob\l\'oj and Spoida \cite{OblSp} that if
their Assumption $\circledast$ fails, then bound {(\ref
{eqlookbackotbound})} is not necessarily optimal, and Theorem~\ref{thm-nmarginal}(ii) fails.
\end{Remark}

%
%re3.7 #&#
\begin{Remark}\label{rksimplephi}
{\rm In the case $\phi=\sum_{j=1}^J {\mathbf1}_{[m_j,\infty)}$ for some
$m_j\ge X_0$ with $\zeta^*_j(m_j)<m_j$, the conditions of Theorem~\ref
{thm-nmarginal} are all easily satisfied; that is, both the local
integrability condition above and the requirement $\hat\lambda_i\in
\mathbb{L}^1(\mu_i)$ are immediate for $i=1,\ldots,n$.
%Moreover, the dynamic hedging strategy $\hat H$ is a simple strategy
%with finite number of jumps. Then, then continuous-time local
%martingale $Y^{\hat H}$ in \eqref{XH}, $Y_0=0$, reduces essentially to
%a finite discrete-time local martingale. Since $Y^{\hat H}_{t_n}\ge
%\phi\big(X^*_{t_n}\big)-\sum_{i=1}^n\hat\lambda_i(X_{t_i})$, by the
%robust superhedging condition, we see that $\big(Y^{\hat H}_{t_n}
%\big)^-\in\mathbb{L}^1(\P)$ for all $\P\in\Pc^\mu$. Then, it follows
%from Lemma~\ref{lemlocalmart} reported in Appendix that $Y^{\hat H}$
%is $\P$-martingale for all $\P\in\Pc^\mu$. Hence, the condition $\hat H
%\in\Hc^\mu$ is also satisfied.
}
\end{Remark}

%s4 #&#
\section{Proofs of the main results}\label{sectpathwise}

%s4.1 #&#
\subsection{Proof of the pathwise inequality of Proposition \texorpdfstring{\protect\ref{lemTrajectorialInequalityOrderedCase}}{3.1}}\label{secttrajectorial}
Our objective is to prove by induction the following trajectorial inequality:
\begin{eqnarray*}
{\mathbf1}_{ \{ {\omega^*_{t_n}\ge m}  \} } &\le& \Upsilon_n (\omega, m ,\zeta)
\\
&:=& \sum
_{i=1}^n \biggl( \frac{(\omega_{t_i}-\zeta_i)^+}{
m-\zeta_i} +{
\mathbf1}_{ \{ {\omega^*_{t_{i-1}<m\le\omega^*_{t_i}}}  \} } \frac{(m-\omega_{t_i})}{
m-\zeta_i} \biggr)
\\
&&{} +\sum_{i=1}^{n-1} \biggl(
\frac{(\omega_{t_i}-\zeta_{i+1})^+}{
m-\zeta_{i+1}} +{\mathbf1}_{ \{ {m\le\omega^*_{t_i},\zeta_{i+1}\le\omega
^*_{t_i}}  \} } \frac{(\omega_{t_{i+1}}-\omega_{t_i})}{
m-\zeta_{i+1}} \biggr),
\end{eqnarray*}
which is immediately seen to imply the required inequality.

For simplicity we omit the arguments $(\omega, m ,\zeta)$ for
$\Upsilon_n$ below.
First, in the case $n=1$, the required inequality is the same as that
of Lemma 2.1 of Brown, Hobson and Rogers \cite{brownhobsonrogers}:
%
%e4.1 #&#
%e4.2 #&#
\begin{eqnarray} \label{eqpwproofindstart}
\Upsilon_1 &=& \frac{(\omega_{t_1}-\zeta_1)^+
+{\mathbf1}_{ \{ {\omega^*_{t_0}<m\le\omega^*_{t_1}}  \} }
(m-\omega_{t_1})}{
m-\zeta_1}\nonumber
\\
&\geq& \frac{\omega_{t_1}-\zeta_1
+m-\omega_{t_1}}{
m-\zeta_1} {
\mathbf1}_{ \{ {m\le\omega^*_{t_1}}  \} }
\\
&\ge& {\mathbf1}_{ \{ {m\le\omega^*_{t_1}}  \} }.\nonumber
\end{eqnarray}
We next assume that $\Upsilon_{n-1}\ge{\mathbf1}_{ \{ {\omega
^*_{t_{n-1}}\ge m}  \} }$ for some $n\ge2$, and show that
$\Upsilon_n\ge{\mathbf1}_{ \{ {\omega^*_{t_n}\ge m}  \} }$.
We consider two cases:

\begin{longlist}[\textit{Case} 2:]
\item[\textit{Case} 1: $\omega^*_{t_{n-1}}\ge m$.] Then\vspace*{1pt} $\omega
^*_{t_n}\ge m$, and it follows from the induction hypothesis that
$1={\mathbf1}_{ \{ {\omega^*_{t_n}\ge m}  \} } ={\mathbf1}_{
\{ {\omega^*_{t_{n-1}}\ge m}  \} }\le\Upsilon_{n-1}$. In order
to see that $\Upsilon_{n-1}\le\Upsilon_n$, we compute directly that,
in the present case,
%
%e4.3 #&#
\begin{eqnarray}
\label{eqrebalancingcalendarspreads} \Upsilon_n-\Upsilon_{n-1} &=&
\frac{\omega_{t_n}-\zeta_n}{m-\zeta_n} ({\mathbf1}_{ \{ {\omega_{t_n}\ge\zeta_n}  \} } -{\mathbf1}_{ \{ {\omega_{t_{n-1}}\ge\zeta_{n}}  \} } ) \ge 0.
\end{eqnarray}

\item[\textit{Case} 2: $\omega^*_{t_{n-1}}<m$.] As $(\omega^*_t)$ is\vspace*{2pt}
nondecreasing, it follows that $\omega^*_{t_i}<m$ for all $i\le n-1$.
With a direct computation we obtain
\begin{eqnarray}
\Upsilon_n = \Upsilon^0_n+\frac{(\omega_{t_n}-\zeta_n)^+}{
m-\zeta_n}
+{\mathbf1}_{ \{ {m\le\omega^*_{t_n}}  \} }\frac{m-\omega_{t_n}}{
m-\zeta_n}\nonumber
\\
\eqntext{\displaystyle \mbox{where } \Upsilon^0_n:=
\sum_{i=1}^{n-1} \biggl(\frac{(\omega_{t_i}-\zeta_i)^+}{
m-\zeta_i} -
\frac{(\omega_{t_i}-\zeta_{i+1})^+}{
m-\zeta_{i+1}} \biggr).}
\end{eqnarray}
Since $m>\omega^*_{t_i}\ge\omega_{t_i}$ for $i\le n-1$, the
functions $\zeta\longmapsto(\omega_{t_i}-\zeta)^+/(m-\zeta)$ are
nonincreasing. This implies that $\Upsilon^0_n\ge0$ by the fact that
$\zeta_i\le\zeta_{i+1}$ for all $i\le n$. Then
%
%e4.4 #&#
\begin{eqnarray}\label{eqforwardtransaction}
\Upsilon_n & \ge & \frac{(\omega_{t_n}-\zeta_n)^+
+{\mathbf1}_{ \{ {m\le\omega^*_{t_n}}  \} }(m-\omega_{t_n})}{
m-\zeta_n}
\nonumber
\\
&\ge& \frac{(\omega_{t_n}-\zeta_n)^+
+m-\omega_{t_n}}{
m-\zeta_n} {\mathbf1}_{ \{ {m\le\omega^*_{t_n}}  \} }
\ge \frac{\omega_{t_n}-\zeta_n
+m-\omega_{t_n}}{
m-\zeta_n} {\mathbf1}_{ \{ {m\le\omega^*_{t_n}}  \} }
\\
& = & {\mathbf1}_{ \{ {m\le\omega^*_{t_n}}  \} }.
\nonumber
\end{eqnarray}
\end{longlist}\vspace*{-15pt}
\hspace*{342pt}\qed%
%%\end{pf}

%s4.2 #&#
\subsection{The iterated Az\'ema--Yor-type embedding of Ob{\l}{\'o}j and Spoida \texorpdfstring{\cite{OblSp}}{[41]}}\label{55555222222222}

Before we proceed to the proof of Theorems~\ref{thm-nmarginal-pricing}
and~\ref{thm-nmarginal}, we recall the iterated Az\'ema--Yor-type
embedding of Ob{\l}{\'o}j and Spoida \cite{OblSp}. This embedding
will allow us to identify the extremal model in the context of these theorems.

Under their Assumption $\circledast$, Ob\l\'oj and Spoida
\cite{OblSp} extend the Az\'ema--Yor embedding for $\mu_1,\ldots,\mu
_n$ by introducing the stopping times based on some functions $\eta
_1,\ldots,\eta_n$,
%
%e4.5 #&#
\begin{equation}
\label{eqiteratedAY} \tau_0=0\quad\mbox{and}\quad  \tau_i:= \inf
\bigl\{ t \geq\tau_{i-1}\dvtx  X_t \leq\eta_{i}
\bigl(X^*_t\bigr) \bigr\}, \qquad i=1,\ldots,n.
\end{equation}
Theorem 2.6 therein asserts that for $\eta$ obtained from an iterative
optimization problem (these functions are called $\xi_1,\ldots,\xi_n$
in \cite{OblSp}), we have $X_{\tau_i}\sim_{\mathbb{P}_0}\mu_i$,
and $(X_{t\land\tau_n})$ is uniformly integrable. Consider a time
change of $X$.
%
%e4.6 #&#
\begin{equation}
\label{eqchangeoftimeAY} Z_t:=X_{\tau_i\land (\tau_{i-1}\lor ({t-t_{i-1}})/({t_{i} - t
}) )}\qquad\mbox{for }
t_{i-1}<t\leq t_i, i=1,\ldots,n
\end{equation}
with $Z_0=X_0$, and observe that $Z$ is a continuous, uniformly
integrable martingale on $[0,t_n]$ with $Z_{t_i}=X_{\tau_i}\sim
_{\mathbb{P}_0}\mu_i$.
As a consequence, the distribution of $Z$, $\mathbb{P}^{\max}:=
\mathbb{P}_0\circ(Z)^{-1}$, is an element of ${\mathcal P}^\mu$.

We shall argue in Appendix~\ref{appproofsiAYembedding} that
%
%e4.7 #&#
\begin{equation}
\qquad\zeta^*_i(m) = \min_{j\geq i} \eta_j(m)
\qquad\mbox{for all } i \leq n, X_0< m <\inf\bigl\{y\dvtx
c_n(y)=0\bigr\}. \label{equniqueminimizers}
\end{equation}
Then note that
%
%e4.8 #&#
\begin{equation}
C(m) = K_n(m) \label{eqidentificationC}
\end{equation}
for the continuous and nonincreasing function $K_n$ defined by Ob\l\'
oj and Spoida \cite{OblSp}; see Lemma 2.14 therein.

Optimality of these iterated Az\'ema--Yor-type embeddings will follow
from the fact that they attain a.s. equality in the pathwise
inequality {(\ref
{eqthmTrajectorialInequalityOrderedCase1})}. We state this in a
greater generality:

%
%le4.1 #&#
\begin{Lemma}[(Pathwise equality)]\label{propPathwiseEquality}
Let $(\eta_i)_{1\le i\le n}$ be nondecreasing, right-continuous
functions, and $\zeta_i:=\min_{j\ge i}\eta_j$. Let $(\tau_i)_{1\le
i\le n}$ be the corresponding stopping times defined by {(\ref
{eqiteratedAY})}, and $Z$ the process defined by {(\ref
{eqchangeoftimeAY})}. Assume that $(X_{t\land\tau_n}\dvtx t\geq0)$ is
$\mathbb{P}_0$-uniformly integrable. Then, for any $m > Z_0$ with
$\eta_n(m) < m$, $Z$ achieves equality in {(\ref
{eqthmTrajectorialInequalityOrderedCase1})},
%
%e4.9 #&#
\begin{eqnarray}
\label{eqrelationzetaxipathwiseequalityresult} {\mathbf1}_{ \{ {Z^*_{t_n}(\omega) \geq m }  \} } &=& \Upsilon_n\bigl(Z(
\omega),m,\zeta(m)\bigr)\qquad\forall\omega\in\Omega_{X_0}.
\end{eqnarray}
\end{Lemma}

\begin{pf}
See Appendix~\ref{appproofsiAYembedding}.
\end{pf}

%s4.3 #&#
\subsection{Proof of Theorem \texorpdfstring{\protect\ref{thm-nmarginal-pricing}}{3.3}\textup{(ii)}}\label{secproofdirect}
Assume that $(\mu_i)$ satisfy Assumption $\circledast$ of
Ob\l\'oj and Spoida \cite{OblSp}. Recall that {(\ref
{equniqueminimizers})} then holds. First, suppose $\zeta_n^*(m)<m$
$d \phi(m)$-a.e. Then it follows directly from Lemma~\ref
{propPathwiseEquality} and the proof of Proposition~\ref
{propmartineq} that we have equality in {(\ref
{eqlookbackotbound})} which is attained by $\mathbb{P}^{\max}:=
\mathbb{P}_0\circ(Z)^{-1}\in{\mathcal P}^\mu$, as defined in
Section~\ref{55555222222222}.
Finally, if $d\phi(m)$ charges the set $\{m\dvtx  \zeta^*_n(m)=m\}$, the
following reasoning applies.
It is seen directly that Assumption $\circledast$ excludes
that $c_n \equiv c_{n-1}$ on an open interval inside the support of
$\mu_n$.
Then for every $\delta>0$ and $m \in( X_0, \inf\{y\dvtx  c_n(y)=0\} ]$
such that $\zeta^*_n(m) = m$ there exists $m' \in(m-\delta,m)$ such
that $\zeta^*_n(m') < m'$.
Consequently, for $\varepsilon> 0$ there exists a right-continuous
nondecreasing $\phi^{\varepsilon}$ such that $0 \leq\phi^{\varepsilon}
- \phi< \varepsilon$, $\phi^\varepsilon(X_0)=\phi(X_0)$ and $\zeta
_n^*(m)<m$ $d\phi^{\varepsilon}$-a.s. Note that clearly $\E^{\mathbb
{P}^{\max}} [ \phi(Z^*_{t_n})  ]$ is finite if and only if
$\E^{\mathbb{P}^{\max}} [ \phi^\varepsilon(Z^*_{t_n})  ]$
is. Let $\phi^{\varepsilon,-1}$ and $\phi^{-1}$ denote the
right-continuous inverses of $\phi^\varepsilon$ and $\phi$,
respectively. Finally, recall from {(\ref{eqidentificationC})}
above that $C$ is nonincreasing and continuous.
Hence, applying the previous case to $\phi^\varepsilon$, we have
\begin{eqnarray*}
\E^{\mathbb{P}^{\max}} \bigl[ \phi\bigl(Z^*_{t_n}\bigr) \bigr] &=& \lim
_{\varepsilon\to0} \E^{\mathbb{P}^{\max}} \bigl[ \phi ^{\varepsilon}
\bigl(Z^*_{t_n}\bigr) \bigr]
\\
& =& \lim_{\varepsilon\to0} \biggl\{
\phi(X_0) + \int_{(X_0,\infty)} C(m) \,d
\phi^{\varepsilon}(m) \biggr\}
\\
& =& \lim_{\varepsilon\to0} \biggl\{ \phi(X_0) + \int
_{\phi
(X_0)}^\infty C\bigl(\phi^{\varepsilon,-1}(x)\bigr) \,d x
\biggr\}
\\
&=& \phi(X_0) + \int_{(X_0,\infty)} C(m) \,d \phi(m),
\end{eqnarray*}
where we used monotone convergence since $C(\phi^{\varepsilon
,-1}(x))\geq C(\phi^{-1}(x))$.
%%%\end{pf}

%s4.4 #&#
\subsection{Proof of Theorem \texorpdfstring{\protect\ref{thm-nmarginal}}{3.5}}\label{secproofdirect1}

Note that by Lemma~\ref{27565734534563}, $\zeta^*_i(m)\geq\ub^{-1}(m)\to
\infty$ as $m\to\infty$, and hence the integral defining $\hat
\lambda_i$ in {(\ref{optimalsemistatic})} is over a bounded
interval. It is therefore well defined by the assumed local
integrability. The same argument, applied to representation {(\ref
{eqHdecomposition})}, shows that $\hat H$ is well defined.
%Finally, the second integral is constant on each interval
%$[t_{i-1},t_i)$, and therefore does not raise any problem for the
%stochastic integral with respect to $X$. As for the first integral, we
%again compute by integration by parts that:
% \b*
% \int_{t_{i-1}}^{t_i}
% \int_{\omega^*_{t_{i-1}}}^{\omega^*_t}
% \frac{d\phi(m)}{m-\zeta^*_i(m)}\,dX_t
% &=&
% \int_{\omega^*_{t_{i-1}}}^{\omega^*_{t_i}}
% \frac{(X_{t_i}-m)\,d\phi(m)}{m-\zeta^*_i(m)},
% \e*
%which is also well-defined by our assumed local integrability
%condition.

The superhedging inequality {(\ref{hedgingstrategy})} then follows
from the trajectorial inequality of Proposition~\ref
{lemTrajectorialInequalityOrderedCase} together with {(\ref
{indicatortophi})}.
If $\hat\lambda\in\Lambda_n^\mu$ and $\hat H\in{\mathcal H}^\mu
$, then {(\ref{hedgingstrategy})} instantly implies the bound
$U^\mu_n (\phi(X^*_T) )\leq\phi(X_0)+ \mu(\hat{\lambda
})$. Finally, by Fubini, $\mu(\hat{\lambda}) = \int_{(X_0,\infty
)}C(m) \,d \phi(m)<\infty$, thus establishing claim (i) of the theorem.

Note that the integrability conditions of Theorem~\ref{thm-nmarginal}
imply, in particular, that $\zeta_n^*(m)<m$ $d \phi(m)$-a.e. If $\mu
$ satisfy Assumption $\circledast$ of Ob\l\'oj and
Spoida \cite{OblSp}, then as above, it follows directly from Lemma
\ref{propPathwiseEquality} by integrating the pathwise equality
against $d \phi$, that there is $\mathbb{P}^{\max}$-a.s. equality
in~{(\ref{hedgingstrategy})}. As a consequence the equality
$U_n^\mu(\phi(X^*_T))=\phi(X_0)+ \mu(\hat\lambda)$ holds. This
establishes claim (ii) of the theorem.

It remains to argue the admissibility of $\hat\lambda$ and $\hat H$.
Observe that for $i=n$, {(\ref{eqfinitelambdacondition})} is
simply the required property $\int\hat\lambda_n \,d\mu_n<\infty$.
Note also that the inner integral in~{(\ref
{eqfinitelambdacondition})} is a convex function of $x$ so that by
{(\ref{ci-nondecreasing})}, we may replace $\mu_i$ with $\mu
_{i-1}$ in~{(\ref{eqfinitelambdacondition})} and the double
integral remains finite. The equivalence of $\hat\lambda\in\Lambda
_n^\mu$ and {(\ref{eqfinitelambdacondition})} now follows from
the definition of $\hat\lambda$.

Finally, we show $\hat H\in{\mathcal H}^\mu$. It is immediate that
$\hat H \in\hat{\mathbb{H}}^2$.
It remains to prove that $\int_0^\cdot\hat H_s \,dX_s$ is a $\mathbb
{P}$-supermartingale for any $\mathbb{P}\in{\mathcal P}^\mu$. Let us
fix one such $\mathbb{P}$ and recall that $X$ is a $\mathbb
{P}$-continuous martingale. For $t\in[t_{i},t_{i+1})$ we have, by {(\ref{eqHdecomposition})},
%
%e4.10 #&#
\begin{equation}
\label{eqintegralofHexplicit} \qquad\int_{t_{i}}^t \hat H_s
\,dX_s = -\int_{t_{i}}^t
\bigl(f_i\bigl(X^*_s\bigr)-f_i
\bigl(X^*_{t_{i}}\bigr)\bigr) \,dX_s -g_i
\bigl(X_{t_i},X^*_{t_i}\bigr) (X_t-X_{t_{i}}),
\end{equation}
where $f_i(x):= \int_{(X_0,x]}\frac{d\phi(m)}{m-\zeta^*_{i+1}(m)}$
and $g_i(x,y) = \int_{(X_0,y]}{\mathbf1}_{ \{ {x\geq\zeta
^*_{i+1}(m)}  \} }\frac{d\phi(m)}{m-\zeta^*_{i+1}(m)}$. In the
first stochastic integral we recognize an Az\'ema--Yor process; see
Carraro, El Karoui and Ob\l\'oj \cite{CKO12}. We recall that
$M^{F_i}(X)_t := F_i (X^*_t) - f_i(X^*_t)(X^*_t-X_t)$,
where $F_i(x)=\int_{X_0}^x f_i(u)\,du=\int_{(X_0,x]}\frac{(x-m)\,d\phi
(m)}{m-\zeta^*_{i+1}(m)}$, $x\geq X_0$, satisfies $M^{F_i}(X)_t = \int_0^t f_i(X^*_s)\,dX_s$ and is a local martingale.

We can extend $F_i$ to the real line putting $F_i(x)=0$ for $x<X_0$.
This preserves convexity, and we observe that we hence have
$M^{F_i}(X)_t \leq F_i(X_t)$. Also, for $t\leq t_{i+1}$,
\begin{eqnarray*}
\E^\mathbb{P}\bigl[ F_i(X_t)\bigr]&\leq&
\E^\mathbb{P}\bigl[F_i(X_{t_{i+1}})\bigr]
\\
 &=& \int
_\mathbb{R}\int_{X_0}^\infty
\frac{(x-m)^+\,d\phi(m)}{m-\zeta
^*_{i+1}(m)}\mu_{i+1}(dx)
\\
&\leq&\int_\mathbb{R}\int_{X_0}^\infty
\frac{(x-\zeta^*_{i+1}(m))^+\,d\phi(m)}{m-\zeta^*_{i+1}(m)}\mu _{i+1}(dx)<\infty
\end{eqnarray*}
by condition {(\ref{eqfinitelambdacondition})}. It follows that
$M^{F_i}(X)_t^+$ is $\mathbb{P}$-integrable.
% and $\E^\P[M^{F_i}(X)_t]$ is well defined.

Let $(\tau_k)$ be a localizing sequence for $M^{F_i}(X)$ and $0\leq
s\leq t\leq t_n$. A simple monotone convergence argument shows that $\E
^\mathbb{P}[F_i(X^*_{t\wedge\tau_k})\mid{\mathcal F}_s]$ converge to
$\E^\mathbb{P}[F_i(X^*_{t})\mid {\mathcal F}_s]$, as $k\to\infty$, and
$M^{F_i}(X)_{s\land\tau_k}$ converge to $M^{F_i}(X)_s$ by continuity.
It follows that the following limit is well defined:
\[
\lim_{k\to\infty} \E^\mathbb{P}\bigl[f_i
\bigl(X^*_{t\wedge\tau
_k}\bigr) \bigl(X^*_{t\wedge\tau_k}-X_{t\wedge\tau_k}\bigr)
\bigl\llvert {\mathcal F}_s\bigr]\geq\E ^\mathbb{P}
\bigl[f_i\bigl(X^*_t\bigr) \bigl(X^*_t-X_t
\bigr)\bigr\rrvert {\mathcal F}_s\bigr]\qquad\mbox{a.s.},
\]
and the inequality follows from Fatou's lemma since $f_i\geq0$ and
$X^*\geq X$. Combining these we obtain
\begin{eqnarray*}
\E^\mathbb{P}\bigl[M^{F_i}(X)_t\mid  {
\mathcal F}_s\bigr]&\geq&\lim_{k\to\infty}
\E^\mathbb{P}\bigl[F_i\bigl(X^*_{t\wedge\tau_k}
\bigr)-f_i\bigl(X^*_{t\wedge\tau
_k}\bigr) \bigl(X^*_{t\wedge\tau_k}-X_{t\wedge\tau_k}
\bigr)\mid  {\mathcal F}_s\bigr]
\\
&=&M^{F_i}(X)_s,
\end{eqnarray*}
where the LHS is well defined since $M^{F_i}(X)_t^+$ is $\mathbb
{P}$-integrable. In particular $\E^\mathbb{P}[M^{F_i}(X)_t]\geq
M^{F_i}(X)_0=0$, and hence $\E^\mathbb{P}[\llvert  M^{F_i}(X)_t\rrvert  ]<\infty$,
and $M^{F_i}(X)$ is a submartingale. Finally, $-M^{F_i}(X)$ is a
supermartingale.

We can compute explicitly
\[
\int_{t_{i}}^t \bigl(f_i
\bigl(X^*_s\bigr)-f_i\bigl(X^*_{t_{i}}\bigr)
\bigr) \,dX_s = M^{F_i}(X)_t -
M^{F_i}(X)_{t_i} - f_i\bigl(X^*_{t_i}
\bigr) (X_t-X_{t_i}).
\]
Combining with {(\ref{eqintegralofHexplicit})} we conclude that
\begin{eqnarray*}
\int_0^{t_n} \hat H_s
\,dX_s &=& - \sum_{i=0}^{n-1}
\bigl(M^{F_i}(X)_{t_{i+1}}-M^{F_i}(X)_{t_{i}}
\bigr)
\\
&&{} + \sum_{i=0}^{n-1}
\bigl(f_i\bigl(X^*_{t_i}\bigr)-g_i
\bigl(X_{t_i},X^*_{t_i}\bigr) \bigr) (X_{t_{i+1}}-X_{t_i}),
\end{eqnarray*}
and we note that by the above, the first sum is integrable under
$\mathbb{P}$. By the superhedging property {(\ref
{hedgingstrategy})}, we have
\begin{eqnarray*}
&& \sum_{i=0}^{n-1} \bigl(f_i
\bigl(X^*_{t_i}\bigr)-g_i\bigl(X_{t_i},X^*_{t_i}
\bigr) \bigr) (X_{t_{i+1}}-X_{t_i})
\\
&&\qquad \geq\phi\bigl(X^*_{t_n}
\bigr)-\hat\lambda(X_\tbf)+\sum_{i=0}^{n-1}
\bigl(M^{F_i}(X)_{t_{i+1}}-M^{F_i}(X)_{t_{i}}
\bigr),
\end{eqnarray*}
and the RHS is a $\mathbb{P}$-integrable r.v. by the above, the
assumption that $\hat\lambda\in\Lambda_n^\mu$ and, as argued
above, its implication that the bound in {(\ref
{eqlookbackotbound})} is finite. It follows from Lemma~\ref
{lemlocalmart} that the simple discrete trading component of $\hat H$
defines a $\mathbb{P}$-martingale; that is, $
(f_i(X^*_{t_i})-g_i(X_{t_i},X^*_{t_i}) )(X_{t}-X_{t_i})$, $t\in
[t_i,t_{i+1})$ is a $\mathbb{P}$-martingale. The above was carried
under an arbitrary $\mathbb{P}\in{\mathcal P}^\mu$ and hence $\hat
H\in{\mathcal H}^\mu$ as required. This completes the proof of the theorem.
%%%%\end{pf}

%s5 #&#
\section{The stochastic control approach}
\label{sectstochasticcontrol}

We now present the methodology which led us to identify the remarkable
pathwise inequality of Proposition~\ref
{lemTrajectorialInequalityOrderedCase}, and to deduce the value
{(\ref{eqglobaloptimizationproblem})} as the cheapest
semi-static superhedging cost.

For technical reasons and clarity of presentation, we impose the
following additional conditions on the nature of the lookback payoff
and the marginal constraints:
%
%e5.1 #&#
%e5.2 #&#
\begin{eqnarray}\label{hyp2ndapproach}
\qquad &&\textstyle  \phi\in C^1 \mbox{ bounded
nondecreasing, }\phi_{\mid(-\infty,X_0]} \equiv0, \int_{X_0}^\infty
C(m)\,d\phi(m) < \infty,
\nonumber\\[-9pt]\\[-9pt]\nonumber
&& \hat\lambda \mbox{ bounded and } \zeta_i^*\mbox{ continuous
increasing}, i=1,\ldots,n,
\end{eqnarray}
where $\hat\lambda$ is the optimal static hedging payoff function of
{(\ref{optimalsemistatic})}.
The assumption that $\phi$ is constant on $(-\infty,X_0]$ is no loss
generality since these values of $\phi$ are irrelevant for the payoff
$\phi(X^*_{T})$. Clearly adding a constant to $\phi$ does not change
the problem, and hence we take $\phi(X_0)=0$ for convenience.

%
%re5.1 #&#
\begin{Remark}{\rm
Given the expression of $\hat\lambda$ as difference of call option's
payoffs, the boundedness condition imposed above may seem
inappropriate. However, we observe that when the probability measures
$\mu_i$ have bounded support, the values taken by $\hat\lambda$
outside all supports are irrelevant. Therefore, we may re-define $\hat
\lambda$ as a bounded function.
}
\end{Remark}

Since the candidate optimal static hedging strategy $\hat\lambda$ is
assumed to be bounded, the robust superhedging problem is not changed
by restricting to bounded static strategies. In the present section, we
even seek more simplification, and we also analyze a slightly different
formulation,
%
%e5.3 #&#
\begin{eqnarray}\label{overlineU}
\overline{U}_n^\mu(\xi) &:=& \inf \bigl\{
Y_0\dvtx \exists (\lambda,H)\in\bigl(\mathbb{L}^\infty
\bigr)^n\times{\mathcal H}, \overline{Y}{}^{H,\lambda}_T
\ge\xi,
\nonumber\\[-8pt]\\[-8pt]\nonumber
&&\hspace*{94pt}\mathbb{P}\mbox{-a.s. for all }\mathbb{P}\in{\mathcal P}_S
\bigr\},
\end{eqnarray}
where
${\mathcal P}_S$ is the subset of ${\mathcal P}$, consisting of
probability measures
\begin{eqnarray*}
&&\mathbb{P}:=\mathbb{P}_0\circ\bigl(X^\alpha
\bigr)^{-1}\qquad\mbox{with } X^\alpha:=\int_0^.
\alpha_s\,dB_s\mbox{ for some } \alpha\in
\mathbb{H}^2(\mathbb{P}_0),
\end{eqnarray*}
and where ${\mathcal H}$ is the subset of dynamic trading strategies
$H\in\hat{\mathbb{H}}^2$ such that the corresponding value process
$Y^H$ is a $\mathbb{P}$-supermartingale for all $\mathbb{P}\in
{\mathcal P}_S$. We note that
%
%e5.4 #&#
\begin{eqnarray}\label{eqscchainineq}
\overline{U}_n^\mu(\xi) &\geq&\inf \bigl\{
Y_0\dvtx \exists (\lambda ,H)\in\bigl(\mathbb{L}^\infty
\bigr)^n\times{\mathcal H}^\mu_S,
\overline{Y}{}^{H,\lambda}_T\ge\xi,  \mathbb{P}\mbox{-a.s. } \forall
\mathbb{P}\in{\mathcal P}_S \bigr\}
\nonumber\\[-8pt]\\[-8pt]\nonumber
&\geq&\sup_{\mathbb{P}\in{\mathcal P}_S^\mu}
\E^\mathbb {P}[\xi],
\end{eqnarray}
where ${\mathcal P}_S^\mu:=\{\mathbb{P}\in{\mathcal P}_S\dvtx X_{t_i}\sim
_\mathbb{P}\mu_i,i=1,\ldots,n\}$ and ${\mathcal H}^\mu_S$ is
obtained by relaxing the supermartingale requirement in ${\mathcal H}$
to $\mathbb{P}\in{\mathcal P}_S^\mu$. The first inequality then
follows by relaxation of conditions on $(\lambda,H)$, and the second
one is the usual majorization of pricing by hedging; see {(\ref
{eqintrophineq})}.
We note that, given the definition of $U^\mu_n$ in {(\ref{V})}
and Theorem~\ref{thm-nmarginal}, the middle term may appear more
natural for the superhedging price. However, under some simplifying
technical assumptions and for $\xi=\phi(X^*_T)$, we will show that in
fact we have equalities throughout {(\ref{eqscchainineq})}.

%s5.1 #&#
\subsection{Dual formulation of the robust superhedging problem}

The first step for the present approach is the observation that
\begin{eqnarray*}
\overline{U}_n^\mu(\xi) &=& \inf_{\lambda\in(\mathbb{L}^\infty)^n}
\inf \bigl\{ Y_0\dvtx \exists H\in{\mathcal H}, Y^H_T
\ge\xi-\lambda(X_{\tbf})+\mu(\lambda),
\\
&&\hspace*{133pt} \mathbb{P}\mbox{-a.s. for all }
\mathbb{P}\in{\mathcal P}_S \bigr\}.
\end{eqnarray*}
We shall continue analyzing the RHS of the last inequality,
observing that for each fixed $\lambda\in(\mathbb{L}^\infty)^n$, we
are reduced to the robust superhedging problem of the derivative $\xi
-\lambda(X_{\tbf})+\mu(\lambda)$. A dual formulation of this
problem is derived in Theorem 2.1 of \cite{ght} under some uniform
continuity assumptions. More recently, Neufeld and Nutz \cite
{Neufeld-Nutz} relaxed the uniform continuity condition, allowing for a
larger class of random variables including measurable ones. The
following direct application of \cite{Possamai-Royer-Touzi} is better
suited to our context:

%
%pr5.2 #&#
\begin{Proposition}\label{propduality}
For a bounded random variable $\xi$, we have
\begin{eqnarray*}
\overline{U}_n^\mu(\xi) &=& \inf_{\lambda\in(\mathbb{L}^\infty)^n}
\sup_{\mathbb{P}\in{\mathcal P}_S} \bigl\{ \mu(\lambda)+\E^\mathbb{P} \bigl[
\xi-\lambda(X_{\tbf
}) \bigr] \bigr\}.
\end{eqnarray*}
\end{Proposition}

%s5.2 #&#
\subsection{The one-marginal problem}\label{sectonemarginal}

We start with an essential ingredient, namely a general one-marginal
construction, which allows us to move from $(n-1)$ to $n$ marginals.

For an inherited maximum $M_0\ge X_0$, we introduce the process
\begin{eqnarray*}
&&M_t := M_0\vee X^*_t\qquad\mbox{for } t\ge0.
\end{eqnarray*}
The process $(X,M)$ takes values in the state space $\Dbf:=\{(x,m)\in
\mathbb{R}^2\dvtx x\le m\}$. Our interest in this section is on the
upper bound on the price of the one-marginal ($n=1$) lookback option
defined by the
payoff
%
%e5.5 #&#
\begin{equation}
\label{xiLookback} \xi= g\bigl(X_T,X^*_T\bigr)\qquad\mbox{for
some } g\dvtx \mathbb{R}\times\mathbb{R}\longrightarrow\mathbb{R}.
\end{equation}

{\renewcommand{\theass}{A}
\begin{ass}\label{assA}
Function $g\dvtx \mathbb{R}
\times\mathbb{R}\longrightarrow\mathbb{R}$ is bounded, $C^1$ in
$m$, absolutely continuous in $x$, and $g_{xx}$ exists as a measure.
\end{ass}}

{\renewcommand{\theass}{B}
\begin{ass}\label{assB}
The function $x\longmapsto
\frac{g_m(x,m)}{m-x}$ is nondecreasing.
\end{ass}}%

For a bounded measurable function $\lambda\dvtx \mathbb{R}\longrightarrow
\mathbb{R}$, we denote $g^\lambda
:= g-\lambda$. Similarly to Proposition 3.1 in Galichon, Henry-Labord\`
ere and Touzi \cite{ght}, it follows from the Dambis--Dubins--Schwarz time
change theorem that the model-free upper bound can be converted
into
%
%e5.6 #&#
\begin{eqnarray}\label{U0}
\overline{U}{}^\mu_1(\xi) = \inf
_{\lambda\in\mathbb{L}^\infty} \sup_{\tau\in{\mathcal T}} \bigl\{\mu(\lambda)+J(
\lambda,\tau) \bigr\}
\nonumber\\[-8pt]\\[-8pt]
\eqntext{\displaystyle\mbox{where } J(\lambda,\tau) := \E^{\mathbb{P}_0}
\bigl[g^\lambda\bigl(X_\tau,X^*_\tau\bigr) \bigr],}
\end{eqnarray}
and ${\mathcal T}$ is the collection of all stopping times $\tau$ such that
%
%e5.7 #&#
\begin{equation}
\label{Tinfty+} \{X_{t\wedge\tau},t\ge0\}\qquad\mbox{is a $\mathbb{P}_0$-uniformly
integrable martingale.}
\end{equation}
Then for every fixed multiplier $\lambda\in\mathbb{L}^\infty$, we
are facing the infinite horizon optimal stopping problem
%
%e5.8 #&#
\begin{equation}
\label{ulambda} u^\lambda(x,m) := \sup_{\tau\in{\mathcal T}}
\E^{\mathbb{P}_0}_{x,m} \bigl[g^\lambda(X_\tau,M_\tau)
\bigr],\qquad (x,m)\in\Dbf,
\end{equation}
where $\E^{\mathbb{P}_0}_{x,m}$ denotes the conditional expectation operator
$\E^{\mathbb{P}_0}[\cdot\mid(X_0,M_0)=(x,m)]$. The dynamic programming
equation corresponding to the optimal
stopping problem $u^\lambda$ defined in {(\ref{ulambda})} is
%
%e5.9 #&#
%e5.10 #&#
\begin{eqnarray}
\label{DPE} \min \bigl\{ u-g^\lambda, -u_{xx} \bigr\}&=&0\qquad \mbox{for }(x,m)\in\Dbf,
\nonumber\\[-8pt]\\[-8pt]\nonumber
u_m(m,m)&=&0\qquad\mbox{for }m\in\mathbb{R}.
\end{eqnarray}
It is then natural to introduce a candidate solution for the dynamic
programming equation defined by a free boundary $\{x=\psi(m)\}$, for
some convenient function $\psi$,
%
%e5.11 #&#
%e5.12 #&#
\begin{eqnarray}
v^\psi(x,m) &=& g^\lambda\bigl(x\wedge\psi(m),m\bigr)+\bigl(x-
\psi(m)\bigr)^+g^\lambda_x\bigl(\psi(m),m\bigr) \label{vpsi0}
\\
&=& g^\lambda(x,m) -\int_{\psi(m)}^{x\vee\psi(m)} (x-
\xi)g^\lambda_{xx}(d\xi,m), \qquad x\le m, \label{vpsi}
\end{eqnarray}
where $g_x^\lambda$ denotes the right-derivative of $g^\lambda$ with
respect to $x$, and $g^\lambda_{xx}$ is the corresponding second
derivative in the sense of distributions. The existence of these
derivatives is justified by the restriction of the function $\lambda$
to the set $\hat\Lambda^\mu$ defined in~{(\ref{hatLambdamu})} below.

Here, $v^\psi(\cdot,m)$ coincides with the obstacle $g^\lambda$ before
the exercise boundary $\psi(m)$, and satisfies $v^\psi_{xx}(\cdot,m)=0$
in the continuation region $(\psi(m),m]$. However, the candidate
solution needs to satisfy more conditions. Namely $v^\psi(\cdot,m)$ must
be above the obstacle and concave in $x$ on $(-\infty,m]$, and it needs
to satisfy the Neumann condition in {(\ref{DPE})}.

For this reason, our strategy of proof consists of first restricting
the minimization in {(\ref{U0})} to those multipliers $\lambda$
in the
set
%
%e5.13 #&#
\begin{equation}
\label{hatLambdamu}
\qquad\hat\Lambda := \bigl\{\lambda\in\mathbb{L}^\infty\dvtx
v^\psi\mbox{ concave in }x\mbox{ and } v^\psi\ge
g^\lambda\mbox{ for some }\psi\in\Psi^\lambda \bigr\},
\end{equation}
where the set $\Psi^\lambda$ is defined in {(\ref{Psilambda})}
below so
that our candidate solution $v^\psi$ satisfies the Neumann condition
in {(\ref{DPE})}. Since $v(\cdot,m)=g^\lambda(\cdot,m)$ on $(-\infty
,\psi(m)]$, it follows that $g^\lambda$ is concave on this range,
thus justifying that the second derivative $g^\lambda_{xx}$ is a
well-defined measure for all $\lambda\in\hat\Lambda$. Also, by
Assumption~\ref{assA}, this guarantees that $\lambda''$ is also a well-defined measure.

By formal differentiation of $v^\psi$, the Neumann condition reduces
to the ordinary differential equation
(ODE)
%
%e5.14 #&#
\begin{eqnarray}\label{defgamma} \label{ODE0}
-\psi'g^\lambda_{xx}(\psi,m) =
\gamma(\psi,m)
\nonumber\\[-8pt]\\[-8pt]
\eqntext{\displaystyle\mbox{where } \gamma(x,m) := (m-x)\frac{\partial}{\partial x} \biggl\{
\frac{g_{m}(x,m)}{m-x} \biggr\}}
\end{eqnarray}
exists a.e. in view of Assumption~\ref{assB}. Similarly to Galichon,
Henry-Labord\`ere and Touzi \cite{ght}, we need
for technical reasons to consider this ODE in the relaxed sense.
We then introduce the weak formulation of the ODE
{(\ref{ODE0})},
%
%e5.15 #&#
%e5.16 #&#
\begin{eqnarray}
\label{ODE} \qquad\psi(m)&<&m \qquad\mbox{for all }m \in\mathbb{R},
\nonumber\\[-8pt]\\[-8pt]\nonumber
- \int_{\psi(E)}g^\lambda_{xx}\bigl(\cdot,\psi
^{-1}\bigr) (d\xi) &=& \int_E\gamma(\psi,\cdot) (dm)
\qquad\mbox{for Borel subsets } E\subset\mathbb{R},
\end{eqnarray}
where $\psi$ is chosen in its right-continuous version and is
nondecreasing by the concavity of $g^\lambda$, and the nonnegativity
of $\gamma$ implied by Assumption~\ref{assB}. In fact, we shall restrict to
those $\psi$ which are continuous and (strictly) increasing, so that
the inverse $\psi^{-1}$ is a well-defined continuous increasing
function. This is the reason for the condition on $\zeta^*$ in {(\ref{hyp2ndapproach})}, and this restriction is adopted here for the
sake of technical simplicity.

We introduce the collection of
all relaxed solutions of {(\ref{ODE0})} with the additional
simplifying assumption of continuity
%
%e5.17 #&#
\begin{eqnarray}
\label{Psilambda} \Psi^\lambda &:=& \bigl\{\psi\dvtx \mathbb{R}\to\mathbb{R}
\mbox{ continuous, increasing, and satisfies (\ref{ODE})} \bigr
\}.
\end{eqnarray}
Notice that the ODE {(\ref{defgamma})}, which motivates the relaxation
{(\ref{ODE})}, does not characterize the free boundary $\psi$ as
it is
not complemented by any boundary condition.

%
%re5.3 #&#
\begin{Remark}\label{rem-ODEaffine}
For later use, we observe that {(\ref{ODE})} implies by direct
integration that the function
\[x\longmapsto \lambda(x) -\int_{\psi(X_0)}^x
\int_{X_0}^{\psi^{-1}(y)}\frac{g_m(\psi(\xi),\xi)}{\xi-\psi(\xi
)}\,d\xi \,dy -\int
_{\psi(X_0)}^xg_x\bigl(\xi,
\psi^{-1}(\xi)\bigr)\,d\xi
\]
is affine.
\end{Remark}

%
%pr5.4 #&#
\begin{Proposition}\label{propPeskir}
Let Assumptions \ref{assA} and \ref{assB} hold true. Then
\begin{eqnarray*}
&&u^\lambda\le v^\psi\qquad\mbox{for any } \lambda\in\hat{\Lambda} \mbox{ and } \psi\in\Psi^\lambda.
\end{eqnarray*}
\end{Proposition}

\begin{pf}
 By the definition of $\Psi^\lambda$, the
function $\psi$ that we will
be manipulating has a well-defined continuous increasing inverse. We
proceed in two steps:

\begin{longlist}[\textit{Step} 1.]
\item[\textit{Step} 1.] We first prove that $v^\psi$ is differentiable in $m$
on the diagonal with
%
%e5.18 #&#
\begin{equation}
\label{Neumann} v^\psi_m(m,m)=0\qquad\mbox{for all } m \in
\mathbb{R}.
\end{equation}
Indeed, since $\psi\in\Psi^\lambda$, it follows from Remark~\ref
{rem-ODEaffine} that
\begin{eqnarray*}
\lambda(x) &=& \alpha_0+\alpha_1x+ \int
_{\psi(X_0)}^x \int_{X_0}^{\psi^{-1}(y)}
\frac{g_m(\psi(\xi),\xi)}{\xi-\psi(\xi)}\,d\xi \,dy
\\
&&{}+\int_{\psi(X_0)}^xg_x
\bigl(\xi,\psi^{-1}(\xi)\bigr)\,d\xi
\end{eqnarray*}
for some constants $\alpha_0,\alpha_1$. Plugging this
expression into (\ref{vpsi0}), we see that for $\psi(m)\le x\le m$,
\begin{eqnarray*}
v^\psi(x,m) &=& g\bigl(\psi(m),m\bigr) - \biggl(\alpha_1+
\int_{X_0}^m\frac{g_m(\psi(\xi),\xi)}{\xi-\psi
(\xi)}\,d\xi \biggr)
\bigl(x-\psi(m) \bigr)
\\
&&{} - \biggl(\alpha_0+\alpha_1\psi(m) +\int
_{\psi(X_0)}^{\psi(m)}\int_{X_0}^{\psi^{-1}(y)}
\frac{g_m(\psi(\xi),\xi)}{\xi-\psi(\xi)}\,d\xi \,dy
\\
&&\hspace*{126pt}{} +\int_{\psi(X_0)}^{\psi(m)}g_x
\bigl(\xi,\psi^{-1}(\xi)\bigr)\,d\xi \biggr)
\\
&=& g \bigl(\psi(m),m \bigr)-\alpha_0-\alpha_1x - \int
_{X_0}^m g_m \bigl(\psi(\xi),\xi
\bigr)\frac{x-\psi(\xi)}{
\xi-\psi(\xi)}\,d\xi
\\
&&{} - \int_{\psi(X_0)}^{\psi(m)}g_x
\bigl(\xi,\psi^{-1}(\xi) \bigr)\,d\xi
\\
&=& g\bigl(\psi(X_0),X_0\bigr)+\int_{X_0}^m
g_m\bigl(\psi(\xi),\xi\bigr)\frac{\xi-x}{\xi
-\psi(\xi)} \,d\xi.
\end{eqnarray*}
Since $g$ is $C^1$ in $m$ by Assumption~\ref{assA}, {(\ref{Neumann})}
follows by direct differentiation
with respect to $m$.

\item[\textit{Step} 2.] Let $\tau\in{\mathcal T}$ be arbitrary, and
define the sequence of stopping times
$\tau_n:=\tau\wedge\inf\{t>0\dvtx \llvert  X_t-x\rrvert  >n\}$. Since $v^\psi$ is
concave, it follows from the It\^o--Tanaka formula that
\begin{eqnarray*}
v^\psi(x,m) &\ge& v^\psi(X_{\tau_n},M_{\tau_n})
-\int_0^{\tau_n} v^\psi_{x}(X_t,M_t)\,dX_t
-\int_0^{\tau_n} v^\psi_{m}(X_t,M_t)\,dM_t.
\end{eqnarray*}
Notice that $(M_t-X_t)\,dM_t=0$. Then since $v^\psi_m(m,m)=0$, it
follows that $v^\psi_{m}(X_t,M_t)\,dM_t=v^\psi_{m}(M_t,M_t)\,dM_t=0$,
and therefore
%
%e5.19 #&#
\begin{eqnarray}
\label{ineq-superhedge} v^\psi(x,m) &\ge& v^\psi(X_{\tau_n},M_{\tau_n})
-\int_0^{\tau_n} v^\psi_{x}(X_t,M_t)\,dX_t.
\end{eqnarray}
Taking expectations in the last inequality, we see that
%
%e5.20 #&#
\begin{equation}
\label{ulambdagevpsi} v^\psi(x,m) \ge \E^{\PP_0}_{x,m}
\bigl[v^\psi(X_{\tau_n},M_{\tau_n}) \bigr] \ge
\E^{\PP_0}_{x,m} \bigl[g^\lambda(X_{\tau_n},M_{\tau_n})
\bigr].
\end{equation}
The required result now follows by the dominated convergence, due to
the boundedness of $g$ and $\lambda$ and by the arbitrariness of $\tau
\in{\mathcal T}$.\quad\qed
\end{longlist}\noqed
\end{pf}

%
%re5.5 #&#
\begin{Remark}
Inequality {(\ref{ineq-superhedge})} is the key step in order to
determine the pathwise inequality of Proposition~\ref
{lemTrajectorialInequalityOrderedCase}. Indeed by sending $n\to
\infty$ and taking $\tau=\tau^\psi:=\inf\{t\dvtx X_t\le\psi(M_t)\}$,
we see that
\begin{eqnarray*}
v^\psi(x,m) &\ge& v^\psi(X_\tau,M_\tau)-
\int_0^\tau v^\psi_x(X_t,M_t)\,dX_t
\\
&=& g(M_\tau)-\lambda(X_\tau) -\int_0^\tau
v^\psi_x(X_t,M_t)\,dX_t.
\end{eqnarray*}
This inequality induces the pathwise inequality once we identify the
optimal $\hat\lambda$ and the corresponding free boundary $\hat\psi
$. For the purpose of the pathwise inequality of Proposition~\ref
{lemTrajectorialInequalityOrderedCase}, the optimal superhedging
strategy is identified by using similarly the iterated values functions
$(v_k)_k$ introduced in Section~\ref{sectpreparation} below.
\end{Remark}

%s5.3 #&#
\subsection{Multiple-marginals penalized value function}\label{Multiple-marginals}

We now continue our general methodology and return to the
multiple-marginal problem. Our aim is to prove Theorem~\ref
{thm-nmarginal} for the modified robust superhedging problem $\overline
{U}_n^\mu(\xi)$, and derive the robust superhedging bounds for the
lookback derivative security
%
%e5.21 #&#
\begin{equation}
\qquad \phi\bigl(X^*_T\bigr)\qquad\mbox{given the marginals } X_{t_i}
\sim\mu_i\qquad\mbox{for all }i=1,\ldots,n.
\end{equation}
We recall that the probability measures $\mu_i$ are defined from
market prices which do not admit arbitrage; that is, {(\ref
{ci-nondecreasing})} holds.

Our purpose in this section is to analyze the upper bound on the robust
superhedging cost introduced in Proposition~\ref{propduality},
%
%e5.22 #&#
\begin{equation}
\label{Umun0} \overline{U}{}^\mu_n(\xi) = \inf
_{\lambda\in(\mathbb{L}^\infty)^n} \bigl\{\mu(\lambda)+u^\lambda(X_0,X_0)
\bigr\},
\end{equation}
where
%
%e5.23 #&#
\begin{equation}
u^\lambda(x,m) := \sup_{\mathbb{P}\in{\mathcal P}_S} \E^\mathbb{P}_{x,m}
\bigl[\phi^\lambda(X_\tbf,M_{t_n}) \bigr]\quad\mbox{and}\quad \phi^\lambda:=\phi-\sum_{i=1}^n
\lambda_i.
\end{equation}
Our approach is to introduce the sequence of intermediate optimization problems
%
%e5.24 #&#
\begin{eqnarray}
\label{ulambdak}\qquad  u_n(x,m) =\phi(m)\quad\mbox{and}\quad u_{k-1}(x,m) =
\sup_{\mathbb{P}\in{\mathcal P}_S} \E^\mathbb{P}_{t_{k-1},x,m}
\bigl[u^\lambda_k(X_{t_k},M_{t_k}) \bigr],
\nonumber\\[-10pt]\\[-8pt]
\eqntext{k\le n,}
\end{eqnarray}
where $\E^\mathbb{P}_{t_{k-1},x,m}=\E^\mathbb{P} [\cdot
\mid (X,M)_{t_{k-1}}=(x,m) ]$,
and
%
%e5.25 #&#
\begin{equation}
u^\lambda_k(x,m):=u_k(x,m)-
\lambda_k(x)\qquad\mbox{for } (x,m)\in\Dbf.
\end{equation}
The last iterative sequence of value functions induces $u^\lambda
=u^\lambda_0$. Moreover, by the Dambis--Dubins--Schwarz theorem (see
Proposition 3.1 in \cite{ght}), we may convert the stochastic control
problem in
{(\ref{ulambdak})} into a sequence of optimal stopping problems
%
%e5.26 #&#
\begin{eqnarray}
u_{k-1}(x,m) &=& \sup_{\tau\in{\mathcal T}} \E^{\PP_0}_{x,m}
\bigl[u^\lambda_k(X_\tau,M_\tau) \bigr].
\end{eqnarray}
Then, denoting by
${\mathcal S}_n:=\{\tau=(\tau_1,\ldots,\tau_n)\in{\mathcal
T}\dvtx \tau_1\le
\cdots\le\tau_n\}$, we see that
%
%e5.27 #&#
\begin{eqnarray}\label{Umun}
\overline{U}{}^\mu_n(\xi) = \inf
_{\lambda\in(\mathbb{L}^\infty)^n} \bigl\{ \mu(\lambda)+u^\lambda_0(X_0,X_0)
\bigr\}
\nonumber\\[-8pt]\\[-8pt]
\eqntext{\displaystyle\mbox{with } u^\lambda_0(x,m) := \sup
_{\tau\in{\mathcal S}_n} \E^{\PP_0}_{x,m} \bigl[
\phi^\lambda (X_\tau,M_{\tau
_n} ) \bigr].}
\end{eqnarray}

%s5.4 #&#
\subsection{Preparation for the upper bound}\label{sectpreparation}

The function $u_{k-1}$ corresponds to the optimization problem
considered in Section~\ref{sectonemarginal} with a payoff
$g(x,m)=u_k(x,m)$ depending on the spot and the running maximum.
This was our original motivation for isolating the one-marginal
problem.

To solve the multiple marginals problem, we introduce the iterative
sequence of candidate value functions
%
%e5.28 #&#
\begin{eqnarray}
\qquad v_n(x,m) &:=& \phi(m), \qquad v^\lambda_k(x,m) :=
v_k(x,m)-\lambda_k(x)\quad\mbox{and}
\nonumber
\\
v_{k-1}(x,m) &:=& v^\lambda_k \bigl(x\wedge
\psi_k(m),m\bigr) + \bigl(x-\psi_k(m) \bigr)^+
\partial^+_xv^\lambda_k \bigl(
\psi_k(m),m \bigr) \label{ui}
\\
&=& v^\lambda_k(x,m) -\int_{\psi_k(m)}^{x\vee\psi_k(m)}(x-
\xi)\partial_{xx}v^\lambda _k(d\xi,m),\nonumber
\end{eqnarray}
where $\partial_x^+ v_k^{\lambda}$ and $\partial_{xx} v_k^{\lambda
}$ denote the right-derivative and the measure second derivative,
respectively, of $v_k^{\lambda}$ with respect to $x$ (which are well
defined for $\lambda$ in the subclass $\hat\Lambda_n$ introduced
below), and $\psi=(\psi_1,\ldots,\psi_n)$ with $\psi_i$ defined as
an arbitrary solution of the ordinary differential equation
%
%e5.29 #&#
\begin{eqnarray}\label{odemultiple}
-\psi_k'\partial_{xx}v^\lambda_k(
\psi_k,m) = \gamma_k(\psi_k,m)
\nonumber\\[-8pt]\\[-8pt]
\eqntext{\displaystyle
\mbox{with } \gamma_k(x,m) :=(m-x)\partial_x \biggl\{
\frac{\partial_mv_k(x,m)}{m-x} \biggr\},}
\end{eqnarray}
which stays strictly below the diagonal. Notice that, in contrast to
the one-marginal case, we have dropped here the dependence of $v_k$
in $\psi$ by simply denoting $v_k:=v^\psi_k$ and
$v^\lambda_k:=v^{\psi,\lambda}_k$.

Similarly to the one-marginal case, we introduce the weak formulation
of this ODE,
%
%e5.30 #&#
%e5.31 #&#
\begin{eqnarray}
\label{ODEmiltiple} \qquad\psi_k(m)&<&m\qquad\mbox{for all }m\ge0,
\nonumber\\[-8pt]\\[-8pt]\nonumber
\qquad -\int_{\psi(E)} \partial_{xx}v^\lambda_k
\bigl(\cdot,\psi_k^{-1} \bigr) (d\xi) &=& \int
_E\gamma_k(\psi_k,\cdot) (dm) \qquad
\mbox{for all } E\in{\mathcal B}(\mathbb{R}),
\end{eqnarray}
and we introduce the set
%
%e5.32 #&#
\begin{eqnarray}
\Psi_n^\lambda &:=& \bigl\{\psi\dvtx \mathbb{R}\to
\mathbb{R}^n \mbox{ with continuous increasing entries}
\nonumber\\[-8pt]\\[-8pt]\nonumber
&&\hspace*{94pt}\mbox{$\psi_k$ satisfying (\ref{ODEmiltiple})},k\le n \bigr\}.
\end{eqnarray}
We also follow the one-marginal case by restricting the minimization
in {(\ref{Umun})} to those multipliers $\lambda$ in the set
%
%e5.33 #&#
\begin{eqnarray}
\qquad \hat\Lambda_n &:=& \bigl\{\lambda\in\bigl(\mathbb{L}^\infty
\bigr)^n\dvtx  v_{k-1}\mbox{ concave in $x$ and }v_{k-1}
\ge v^\lambda_k\mbox{ for all }k\le n \bigr\}.
\end{eqnarray}

%
%le5.6 #&#
\begin{Lemma}\label{lem-uv}
Let $\phi$ satisfy {(\ref{hyp2ndapproach})}, $\lambda\in\hat
{\Lambda}_n$ and $\psi\in\Psi^\lambda_n$. Then:

\begin{longlist}[(iii)]
\item[(i)]  for all $i=1,\ldots,n$, the function $v_i$ satisfies
Assumptions \ref{assA} and \ref{assB}, that is, $v_i$ is $C^1$ in $m$, absolutely
continuous in $x$, Lipschitz in
$m$ uniformly in $x$, $\partial_{xx}v_i$ exists as a measure and
$x\longmapsto\partial_mv_i(x,m)/(m-x)$ is nondecreasing;

\item[(ii)]  for all $i=1,\ldots,n$, the function $\partial_mv_i$
is concave in $x$;

\item[(iii)] $u^\lambda(X_0,X_0)\le v_0(X_0,X_0)$.
\end{longlist}
\end{Lemma}

\begin{pf}
 We first prove (i). First $v_n=\phi$
satisfies Assumptions \ref{assA}
and \ref{assB} as it is independent of the $x$-variable, nondecreasing and
$C^1$. For the remaining cases, we proceed by
induction, assuming that $v_{i}$ satisfies Assumptions \ref{assA} and \ref{assB}, for
some $i\le n$,
and we intend to show that $v_{i-1}$ does as well. We first observe
that the following condition is also satisfied by
$v_{i}$:
%
%e5.34 #&#
\begin{equation}
\label{vmconc-recur} v_i(x,m)=\phi(m)\quad\mbox{nondecreasing, or}\quad
\partial_mv_i(m,m)=0,
\end{equation}
where the first alternative holds for $i=n$. $v_{i-1}$ is clearly
$C^1$ in $m$, and by using the ODE {(\ref{odemultiple})} satisfied
by $v_i$,
we directly compute that
%
%e5.35 #&#
\begin{eqnarray}
\label{aboveexpression}
\qquad && \partial_mv_{i-1}(x,m)
\nonumber\\[-2pt]\\[-14pt]\nonumber
&&\qquad  = \cases{
\displaystyle \partial_mv_{i}(x,m), &\quad for $x\in\bigl(-
\infty,\psi_i(m)\bigr]$,
\cr
\displaystyle \partial_mv_{i}
\bigl(\psi_i(m),m\bigr)\frac{m-x}{m-\psi_i(m)}, &\quad for $x\in\bigl[
\psi_i(m),m\bigr]$.}
\end{eqnarray}
Then $v_{i-1}$ inherits the differentiability in $m$, and $x\longmapsto
\partial_mv_{i-1}(x,m)/(m-x)$ is nondecreasing whenever
$x\longmapsto\partial_mv_{i}(x,m)/(m-x)$ is. The remaining properties
follow from the concavity of $v_i$.

We next prove (iii). By the previous step, $v_i$ satisfies
Assumptions \ref{assA} and \ref{assB} for all $i=1,\ldots,n$. Then it follows from
Proposition~\ref{propPeskir} that $u_{n-1}\le v_{n-1}$ for all
$\psi\in\Psi^{\lambda_n}$. Therefore
\begin{eqnarray*}
u_{n-2}(x,m) &\le& \sup_{\tau_{n-1}\in{\mathcal T}} \E^\mathbb{P}_{x,m}
\bigl[v^\lambda_{n-1}\bigl(X_{\tau_{n-1}},X^*_{\tau_{n-1}}
\bigr) \bigr],
\end{eqnarray*}
and we deduce from a second application of Proposition
\ref{propPeskir} that $u_{n-2}\le v_{n-2}$. The required inequality
follows by a backward iteration of this argument.

We finally prove (ii). From {(\ref{aboveexpression})}, we see that
$\partial_mv_{i-1}$ is concave in $x$ on
$ (-\infty,\psi_i(m) )$ and on $(\psi_i(m),m]$. It
remains to verify that $\partial_m v_{i-1}$ is concave at the point
$x=\psi_i(m)$. We directly calculate that
\[
\partial_{xm}v_{i-1} \bigl(\psi_i(m)-,m \bigr)
= \partial_{xm}v_{i} \bigl(\psi_i(m)-,m
\bigr)
\]
and
\[
\partial_{xm}v_{i-1} \bigl(
\psi_i(m)+,m \bigr) = \frac{-\partial_{m}v_{i} ( \psi_i(m),m ) }{
m-\psi_i(m)}.
\]
Then, by the concavity of $\partial_mv_i$ in $x$, together with
{(\ref{vmconc-recur})}, we have
\[
\partial_mv_i \bigl(\psi_i(m),m \bigr) +
\partial_{xm}v_i \bigl(\psi_i(m)+,m \bigr)
\bigl(m-\psi_i(m) \bigr) \ge \partial_mv_i
(m,m ) \ge0,
\]
which implies that $\partial_{xm}v_{i-1} (\psi_i(m)-,m )
\ge\partial_{xm}v_{i-1} (\psi_i(m)+,m )$.
\end{pf}

%
%le5.7 #&#
\begin{Lemma}\label{lem-mu+ulambda}
Let $\phi$ satisfy {(\ref{hyp2ndapproach})}. Then, for
all $\lambda\in\hat{\Lambda}_n$ and $\psi\in\Psi^\lambda_n$, we
have
\begin{eqnarray*}
\mu(\lambda)+u^\lambda(X_0,X_0) &\le& \mu(
\lambda)+v_0(X_0,X_0)
\\
&=& % \phi(X_0)+
\sum_{i=1}^n \int
\bigl[c_i(\xi)-c_0(\xi){\mathbf1}_{\{\xi<\psi_i(X_0)\}} \bigr]
\lambda_i''(d\xi)
\\
&& {}-\int_{\psi_i(X_0)}^{\infty}c_0(
\xi)\partial _{xx}v_i(d\xi,X_0).
\end{eqnarray*}
\end{Lemma}

\begin{pf}
 Denote the LHS$:=\mu(\lambda)+u^\lambda
(X_0,X_0)$. Substituting the expression of the $v_i$'s in inequality
(iii) of Lemma~\ref{lem-uv}, we see that
\begin{eqnarray*}
\mbox{LHS} &\le& % \phi(X_0) +
\sum_{i=1}^n
\mu_i(\lambda_i)-\lambda_i(X_0)
-\int_{\psi_i(X_0)}^{\infty}c_0(\xi)
\partial_{xx}v_i^\lambda(d\xi,X_0)
\\
&\le& % \phi(X_0) +
\sum_{i=1}^n \int
\lambda_i(\xi) (\mu_i-\delta_{X_0}) (d\xi) -
\int_{\psi_i(X_0)}^{\infty}c_0(\xi)
\partial_{xx}v_i^\lambda(d\xi,X_0).
\end{eqnarray*}
Since $\int\xi\mu_i(d\xi)=X_0$, it follows from two integrations by
parts that
\begin{eqnarray*}
\mbox{LHS} &\le& % \phi(X_0) +
\sum_{i=1}^n
\int (c_i-c_0+c_0{\mathbf1}_{[\psi_i(X_0),\infty
)} ) (
\xi)\lambda''_i(d\xi)
\\
&&{} -\int
_{\psi_i(X_0)}^{\infty}c_0(\xi)\partial_{xx}v_i(d
\xi,X_0).
\end{eqnarray*}\upqed
\end{pf}

The following result provides the necessary calculations for the
terms which appear in Lemma~\ref{lem-mu+ulambda}. We denote
%
%e5.36 #&#
\begin{equation}
\overline{\psi}_i := \psi_i\wedge\cdots\wedge
\psi_n\qquad\mbox{for all } i=1,\ldots,n,
\end{equation}
and we set $\overline{\psi}_{n+1}(m):=m$, $m\in\mathbb{R}$.

%
%le5.8 #&#
\begin{Lemma}\label{111122222454}
Let $\phi$ satisfy {(\ref{hyp2ndapproach})}, $\lambda\in\hat
\Lambda_n$, $\psi\in\Psi^\lambda_n$ and $i\le n$. Then for any
$\lambda''_i$-integrable function $\varphi$ we have
\begin{eqnarray}
\int\varphi(\xi)\lambda_i''(d\xi) =
\int \biggl(\frac{\varphi(\overline{\psi}_i(m))}{m-\overline{\psi}_i(m)} -{\mathbf1}_{\{i<n\}}\frac{\varphi(\overline{\psi
}_{i+1}(m))}{m-\overline{\psi}_{i+1}(m)} \biggr){
\mathbf1}_{\{\overline{\psi}_i<\overline{\psi}_{i+1}\}}\,d\phi(m),\nonumber
\\
\eqntext{i=1,\ldots,n.}
\end{eqnarray}
%
%with $\overline{\psi}_{n+1}(m):=m$, $m\in\R$.
\end{Lemma}

\begin{pf}
 See Appendix~\ref{sectAppendix}. \end{pf}

Plugging these calculations into the estimate of Lemma
\ref{lem-mu+ulambda} provides:

%
%le5.9 #&#
\begin{Lemma}\label{lem-mu+un}
Let $\phi\in C^1$ be bounded nondecreasing, $\lambda\in\hat
{\Lambda}_n$ and $\psi\in\Psi^\lambda_n$. Then
\begin{eqnarray*}
&&\mu(\lambda) + u^\lambda(X_0,X_0) \le \mu(
\lambda)+v_0(X_0,X_0) = %\phi(X_0) +
\sum
_{i=1}^n\int_{-\infty}^\infty{
\mathbf1}_{\{\overline{\psi
}_i<\overline{\psi}_{i+1}\}}A_i(m)\,d\phi(m),
\end{eqnarray*}
where
%, denoting $\overline{\psi}_{n+1}(m):=m$, $m\in\R$, and $
%\psi_{1,i}^0:=(\psi_1\wedge\ldots\wedge\psi_i)(X_0)$:
%
\[
A_i := \frac{c_i(\overline{\psi}_i)%-c_0(\overline{\psi}_i) \1_{\{m\le X_0
%\}}\1_{\{\psi_i<\psi_{1,i}^0\}}
}{
m-\overline{\psi}_i} -{\mathbf1}_{\{i<n\}}
\frac{c_i(\overline{\psi}_{i+1}) % -c_0(\overline{\psi}_{i+1})\1_{\{m
%\le X_0\}}\1_{\{\overline{\psi}_{i+1}<\psi_{1,i}^0\}}
}{
m-\overline{\psi}_{i+1}}.
\]
\end{Lemma}

\begin{pf}
We proceed in two steps:

\begin{longlist}[\textit{Step} 1.]
\item[\textit{Step} 1.] We start by reducing the last integral in Lemma~\ref
{lem-mu+ulambda} to an integral with respect to $\lambda_i''$. Notice
that $\partial_{xx}v_i(x,m)={\mathbf1}_{\{i<n\}}{\mathbf1}_{\{x<\psi
_{i+1}(m)\}}\partial_{xx} v_{i+1}^\lambda(x,m)$. Then
\begin{eqnarray*}
&& \int_{\psi_i(X_0)}^\infty c_0(\xi)
\partial_{xx}v_i(d\xi,X_0)
\\
&&\qquad = {
\mathbf1}_{\{i<n\}} \int_{\psi_i(X_0)}^\infty{
\mathbf1}_{\{\xi<\psi_{i+1}(X_0)\}}c_0(\xi )\partial_{xx}v^\lambda_{i+1}(d
\xi,X_0)
\\
&&\qquad = -{\mathbf1}_{\{i<n\}} \int_{\psi_i(X_0)}^\infty{
\mathbf1}_{\{\xi<\psi
_{i+1}(X_0)\}}c_0(\xi)\lambda''_{i+1}(d
\xi)
\\
&&\quad\qquad{} +{\mathbf1}_{\{i<n\}} \int_{\psi_i(X_0)}^\infty{
\mathbf1}_{\{\xi<\psi
_{i+1}(X_0)\}}c_0(\xi)\partial_{xx}v_{i+1}(d
\xi,X_0).
\end{eqnarray*}
By direct iteration, it follows from the fact that $v_n(x,m)=\phi(x)$
is independent of~$x$ that
\begin{eqnarray*}
&& \int_{\psi_i(X_0)}^\infty c_0(\xi)
\partial_{xx}v_i(d\xi,X_0)
\\
&&\qquad = -{
\mathbf1}_{\{i<n\}} \sum_{j=i+1}^n \int
_{\psi_i(X_0)}^\infty{\mathbf1}_{\{\xi<\psi_{i+1}\wedge\psi
_j(X_0)\}}c_0(
\xi)\lambda''_j(d\xi).
\end{eqnarray*}

\item[\textit{Step} 2.] By Lemma~\ref{lem-mu+ulambda}, we have $\mu(\lambda
)+u^\lambda(X_0,X_0)\le A$, where, by the first step,
\begin{eqnarray*}
A &:=& \int \Biggl(\sum_{i=1}^n(c_i-c_0{
\mathbf1}_{(-\infty,\psi_i(X_0)]} )\lambda_i'' +\sum
_{i=1}^n\sum_{j=i+1}^n
c_0{\mathbf1}_{[\psi_i(X_0),\psi
_{i+1}\wedge\cdots\wedge\psi_j(X_0)]} \lambda_j''
\Biggr)
\\
&=& \int \Biggl(\sum_{i=1}^n(c_i-c_0{
\mathbf1}_{(-\infty,\psi_i(X_0)]} )\lambda_i'' +\sum
_{j=1}^n\sum_{i=1}^{j-1}
c_0{\mathbf1}_{[\psi_i(X_0),\psi
_{i+1}\wedge\cdots\wedge\psi_j(X_0)]} \lambda_j''
\Biggr)
\\
&=& \int \Biggl(\sum_{i=1}^n(c_i-c_0{
\mathbf1}_{(-\infty,\psi_i(X_0)]} )\lambda_i'' +\sum
_{j=1}^n c_0{
\mathbf1}_{[\psi_1\wedge\cdots\wedge\psi
_i(X_0),\psi_i(X_0)]} \lambda_j'' \Biggr)
\\
&=& \int\sum_{i=1}^n(c_i-c_0{
\mathbf1}_{(-\infty,\overline\psi
_i(X_0))} )\lambda_i''.
\end{eqnarray*}
Then it follows from Lemma~\ref{111122222454} that $A= \sum_{i=1}^n\int
{\mathbf1}_{\{\overline\psi_i<\overline{\psi}_{i+1}\}} A_i^0\,d\phi$, where
\begin{eqnarray*}
A_i^0 &:=& \frac{c_i(\psi_i)-c_0(\psi_i)
{\mathbf1}_{\{\overline\psi_i<\overline\psi_i(X_0)\}}}{
m-\overline\psi_i}
\\
&&{} -{\mathbf1}_{\{i<n\}}
\frac{c_i(\overline{\psi}_{i+1})
-c_0(\overline{\psi}_{i+1})
{\mathbf1}_{\{\overline{\psi}_{i+1}<\overline\psi_i(X_0)\}}}{
m-\overline{\psi}_{i+1}}.
\end{eqnarray*}
We next observe from the increase of $\overline\psi_i$ that for
$m\geq X_0$, we have $\overline\psi_i(m)\ge\overline\psi_i(X_0)$,
implying that ${\mathbf1}_{\{\overline\psi_i<\overline\psi_i(X_0)\}
}={\mathbf1}_{(-\infty,X_0]}{\mathbf1}_{\{\overline\psi_i<\overline\psi
_i(X_0)\}}$. It also follows that on $\{\overline\psi_i<\overline
{\psi}_{i+1}\}$, we have ${\mathbf1}_{\{\overline{\psi}_{i+1}<\overline
\psi_i(X_0)\}}={\mathbf1}_{(-\infty,X_0]}{\mathbf1}_{\{\overline{\psi
}_{i+1}<\overline\psi_i(X_0)\}}$. The result follows since by {(\ref{hyp2ndapproach})} we have $d \phi(m)\equiv0$ on $(-\infty,X_0]$.\quad\qed
\end{longlist}\noqed
\end{pf}

%s5.5 #&#
\subsection{Proof of \texorpdfstring{$\overline{U}_n^\mu(\xi)=\mu(\hat\lambda)$ under (\protect\ref{hyp2ndapproach}), Assumption $\circledast$ of \cite{OblSp} and Assumptions \protect\ref{assA},~\protect\ref{assB}}%
{Proof of $overline{U}_n^mu(xi)=mu(hat lambda)$ under (5.1), Assumption $circledast$ of [41] and Assumptions A,~B}}

\textit{Step} 1. We first show that $\overline{U}_n^\mu(\xi)\le\mu(\hat
\lambda)$. Given the results of Lemma~\ref{lem-mu+un}, we prove in
this first step that the pointwise minimization % of Lemma~\ref{lem-inf%x<S0} and
in {(\ref{eqglobaloptimizationproblem})} can be achieved by some
vector of Lagrange multipliers
$\hat\lambda=(\hat\lambda,\ldots,\hat\lambda)\in\hat{\Lambda
}_n$, thus
implying that our required upper bound satisfies
%
%e5.37 #&#
\begin{equation}
\label{proofmainresult-ineq1} \overline{U}{}^\mu_n(\xi) \le % \phi(X_0) +
\int
_{X_0}^\infty \sum_{i=1}^n
\biggl( \frac{c_i (\zeta^*_i(m) )}{
m-\zeta^*_i(m)} -\frac{c_i (\zeta^*_{i+1}(m),m )}{
m-\zeta^*_{i+1}(m)} {\mathbf1}_{\{i<n\}} \biggr) \,d
\phi(m).
\end{equation}
In order to define $\hat\lambda$, we take the family of functions
$\hat\psi_i$ given by the boundaries
$\xi_i$ constructed by Ob\l\'oj and Spoida \cite{OblSp}, so that
\begin{eqnarray*}
&&\hat\psi_1:=b_1^{-1}\quad\mbox{and}\quad \overline{
\hat\psi}_i:=\hat\psi_i\wedge\cdots\wedge\hat\psi
_n=\zeta^*_i,\qquad 1<i \leq n.
\end{eqnarray*}
%
%i.e. $\psi^*_i$ is an extension of $\zeta_i^*$ for all $i=1,\ldots,n$.
Recall that $b_1^{-1}$ is the minimizer in {(\ref
{eqone-dimoptimisation})} and is also the right-continuous inverse of
the barycentre function of $\mu_1$; see {(\ref{barycenter})}
below. Also, under Assumption $\circledast$ of Ob\l\'oj and Spoida
\cite{OblSp}, $\hat\psi$ are continuous. Moreover,
direct verification reveals that the functions $\hat\lambda_i$ solve
the system of ODEs {(\ref{odemultiple})}. The required result
follows from our additional assumption in {(\ref
{hyp2ndapproach})} that $\hat\lambda$ is bounded.
%Under our assumption that the support of $\phi'$ is bounded from
%above, it also follows that, up to a linear function, $\lambda_i$ is
%bounded from below. Then
%$\sup_{\P\in\Pc^*}\E^\P\left[\left(\phi(X^*_T)-\lambda(X_{\tbf})
%\right)^+\right]\le\sup_{\P\in\Pc^*}\E^\P\left[\phi(X^*_T)^+\right]<
%\infty$.

\textit{Step} 2. Now we prove that the equality $\overline
{U}_n^\mu(\xi)=\mu(\hat\lambda)$ holds.
In Section~\ref{55555222222222} we recalled the $n$-marginal embedding ${\tau
}_1,\ldots,{\tau}_n$ of Ob\l\'oj and Spoida \cite{OblSp}.
In their Theorem~2.6 they compute the law of $X^{*}_{{\tau}_n}$ as
%
%e5.38 #&#
\begin{equation}
\PP \bigl[ X^*_{{\tau}_n} \geq m \bigr] = K_n(m) = C(m)
%\sum_{i=1}^{n} \left( \frac{c_i( \tilde{\zeta}_i(m))}{ m - \tilde{
%\zeta}_i(m)} - \frac{c_i(\tilde{\zeta}_{i+1}(m))}{m-\tilde{
%\zeta}_{i+1}(m)} \1_{\{i<n\}} \right)
\label{eqdistributionmaximumiAYsecondproof};
\end{equation}
see also {(\ref{equniqueminimizers})} and {(\ref
{eqidentificationC})}.

By definition of $\overline{U}{}^\mu_n(\xi)$ in
{(\ref{overlineU})}, it follows that
%
%e5.39 #&#
\begin{equation}
\overline{U}_n^\mu(\xi) \geq % \inf_{\lambda\in(\mathbb{L}^\infty)^n}
% \Big\{
% \mu(\lambda)
% +\E^{\PP_0}\Big[ \phi\left(X^*_{\tau_n}\right)
% -\sum_{i=1}^n \lambda_i\left(X_{\tau_i}\right)
% \Big]
% \Big\}
% =
\E^{\PP_0} \bigl[ \phi \bigl(X^*_{\tau_n} \bigr) \bigr] % \nonumber\\
\stackrel{{\fontsize{8.36pt}{10pt}{(\ref{eqdistributionmaximumiAYsecondproof})}}} {=} % \phi(X_0) +
\int_{(X_0,\infty)}
C(m) \,d\phi(m). \label{eqlawofmaximum}
\end{equation}
Furthermore, it follows that we have equalities throughout in {\rm
(\ref{eqscchainineq})}.\quad\qed
%The proof is complete by comparing the expression inside the integral
%in \eqref{eqlawofmaximum} with
%\eqref{eqglobaloptimizationproblem} and recalling that $\zeta^{*}_1,
%\dots,\zeta^{*}_n$ was chosen as the minimizer of the latter.

%The remaining claims can be argued as in section
%\ref{secproofdirect}.

%s6 #&#
\section{The Az\'ema--Yor embedding solves the one-marginal problem}
\label{sectAzema-Yor}

In this subsection, we return to the one-marginal context of
Section~\ref{sectonemarginal} with the intention of revisiting the recent
result of Hobson and Klimmek \cite{hobson-klimmek} in this setting.
Our emphasis is again on the efficiency of the stochastic control
approach in the present setting. Therefore, similarly to the previous
section, our assumptions below will be much stronger than what one
could achieve directly with the pathwise approach. The endpoints of the
support of the distribution $\mu$ are denoted
by
%
%e6.1 #&#
\begin{equation}
\label{eqendpointssupport} \ell^\mu:=\sup \bigl\{x\dvtx \mu \bigl([x,\infty) \bigr)=1
\bigr\}\quad\mbox{and}\quad r^\mu:=\inf \bigl\{x\dvtx \mu \bigl((x,\infty) \bigr)=0
\bigr\}.
\end{equation}
We introduce the so-called barycentre function
%
%e6.2 #&#
\begin{equation}
\label{barycenter} b(x) := \frac{\int_{[x,\infty)} y\mu(dy)}{
\mu ([x,\infty) )} {\mathbf1}_{\{x<r^\mu\}} +x {
\mathbf1}_{\{x\ge r^\mu\}}, \qquad x \in\mathbb{R}.
\end{equation}
The solution of Az\'ema and Yor \cite{AzemaYor1,AzemaYor2} to the {\it
Skorokhod embedding problem} is
%
%e6.3 #&#
\begin{eqnarray}
\label{tau*} \hat\tau &:=& \inf \bigl\{t\geq0\dvtx X^*_t\ge
b(X_t) \bigr\} = \inf\bigl\{t\geq0\dvtx  X_t\leq
b^{-1}\bigl(X^*_t\bigr)\bigr\},
\end{eqnarray}
as recalled in Section~\ref{secSEP}.
Plugging $\hat\psi:=b^{-1}$ into the ODE {(\ref{ODE})}, we
obtain the function
%
%e6.4 #&#
\begin{eqnarray}\label{lambda*}
\hat\lambda(x) := \int_{\ell^\mu}^x
\int_{\ell^\mu}^{y} g_m \bigl(\xi,b(\xi)
\bigr) \frac{\mu(d\xi)}{\mu([\xi,\infty))} \,dy +\int_{\ell^\mu}^x
g_x\bigl(\xi,b(\xi)\bigr)\,d\xi ;
\nonumber\\[-8pt]\\[-8pt]
\eqntext{x\in\bigl(-\infty,r^\mu \bigr),}
\end{eqnarray}
whose well-posedness will be guaranteed by the following condition.

{\renewcommand{\theass}{C}
\begin{ass}\label{assC}
$\hat\lambda$ is
well defined and bounded. Moreover, the function $g$ has a measure
second partial derivative with respect to $x$ satisfying
\begin{eqnarray*}
&&g_{xx}(dx,m)-g_{xx}\bigl(dx,b(x)\bigr) \le \gamma
\bigl(x,b(x)\bigr)b'(dx)\qquad\mbox{whenever } b(x)\le m.
\end{eqnarray*}
\end{ass}}%

Similarly to the previous section, we focus on the robust superhedging
problem $\overline{U}{}^\mu_1(\xi)$, as introduced in {(\ref
{overlineU})}, and re-expressed in {(\ref{U0})}.

%
%th6.1 #&#
\begin{Theorem}\label{thmLookback}
Let $\xi=g(X_T,X^*_T)$ for some payoff function $g$
satisfying Assumptions \ref{assA}, \ref{assB} and \ref{assC}. Then the pair $(\hat\lambda,\hat
\tau)$ is a solution of the problem $\overline{U}{}^\mu_1(\xi)$ in~{(\ref{U0})} and
\begin{eqnarray*}
\overline{U}{}^\mu_1(\xi) &=& \mu(\hat\lambda)+ J(\hat
\lambda,\hat\tau) = \E^{\PP_0} \bigl[g \bigl(X_{\hat\tau},X_{\hat\tau}^*
\bigr) \bigr].
\end{eqnarray*}
\end{Theorem}

The remaining part of this section is dedicated to the proof of this
result. Our starting point is the result of Proposition~\ref
{propPeskir} which provides an upper bound for the value function
$\overline{U}{}^\mu_1(\xi)$
for every choice of a multiplier $\lambda\in\hat\Lambda$ and a
corresponding solution $\psi\in\Psi^\lambda$ of the ODE {(\ref{ODE})},
%
%e6.5 #&#
\begin{equation}
\label{upperboundU0} \overline{U}{}^\mu_1(\xi) \le \mu(
\lambda)+v^\psi(X_0,X_0)\qquad\mbox{for all }
\lambda\in\hat\Lambda\mbox{ and } \psi\in\Psi^\lambda.
\end{equation}
Alternatively, for any choice of a nondecreasing function $\psi$
with $\psi(m)<m$ for all $m \in\RR$, we may define a corresponding
multiplier function $\lambda$ by {(\ref{ODE})}, or equivalently by
{(\ref{ODE0})}, in the distribution sense. Then $\psi\in\Psi
^\lambda$.
If in addition $v^\psi$ is concave in $x$ and above the
corresponding\vspace*{1pt} obstacle $g^\lambda$, then $\lambda\in\hat\Lambda$,
and we may conclude by Proposition~\ref{propPeskir} that
$\overline{U}{}^\mu_1(\xi)\le\mu(\lambda) +v^\psi$. The next
result exhibits this bound for the
choice $\psi=b^{-1}$, the right-continuous inverse of the barycentre
function.

%
%pr6.2 #&#
\begin{Proposition}\label{prop-upperbound}
Let $\xi=g(X_T,X^*_T)$ for some payoff function $g$ satisfying
Assumptions \ref{assA}, \ref{assB} and \ref{assC}. Then
\begin{eqnarray*}
\overline{U}{}^\mu_1(\xi) &\le& \mu(\hat\lambda)+J(\hat
\lambda,\hat\tau) = \E^{\PP_0}\bigl[g \bigl(X_{\hat\tau},X_{\hat\tau}^*
\bigr)\bigr].
\end{eqnarray*}
\end{Proposition}

\begin{pf}
 It is immediately checked that
$\hat\psi:=b^{-1}\in\Psi^{\hat\lambda}$. Moreover, by Assumption
\ref{assC} and the subsequent discussion, we see that $\hat\lambda\in\hat
\Lambda$. In view of the previous discussion, the required inequality
follows from Proposition~\ref{propPeskir} once we prove that $v^{\hat
\psi}$ is concave, and that $v^{\hat\psi}\ge g^{\hat\lambda}$.

\begin{longlist}[(2)]
\item[(1)] We first verify that $v^{\hat\psi}$ is concave. By direct
computation using the expression of $\hat\lambda$ in {(\ref{lambda*})}
together with the identity
\begin{eqnarray*}
\frac{b'(dx)}{b(x)-x} &=& \frac{\mu(dx)}{\mu([x,\infty))},
\end{eqnarray*}
we see that
%
%e6.6 #&#
\begin{eqnarray}
g^{\hat\lambda}_{xx}(dx,m) &=& g_{xx}(dx,m)-g_{xx}
\bigl(dx,b(x) \bigr) -\gamma \bigl(x,b(x) \bigr)b'(dx)
\label{glambdaxx}
\end{eqnarray}
in the distribution sense. By Assumption \ref{assC}, it follows that
$x\longmapsto g^{\hat\lambda}(x,m)$ is concave on
$(-\infty,\hat\psi(m)]$. Since $v^{\hat\psi}(\cdot,m)$ is linear
on $[\hat\psi(m),m]$ and $C^1$ across the boundary $\hat\psi$, this
proves that $v^{\hat\psi}$ is
concave.

\item[(2)] We next check that $v^{\hat\psi}\ge g^{\hat\lambda}$. Since
equality holds on $(-\infty,\hat\psi(m)]$, we only compute for
$x\in[\hat\psi(m),m]$ that
\begin{eqnarray*}
\bigl(v^{\hat\psi}-g^{\hat\lambda} \bigr) (x,m) &=& \int
_{\hat\psi(m)}^x \bigl(g^{\hat\lambda}_x
\bigl(\hat\psi(m),m\bigr) -g^{\hat\lambda}_x(\xi,m) \bigr)\,d\xi
\\
&=& -\int_{\hat\psi(m)}^x\int_{\hat\psi(m)}^\xi
g^{\hat\lambda
}_{xx}(dy,m)\,d\xi.
\end{eqnarray*}
By {(\ref{glambdaxx})}, this provides
\begin{eqnarray*}
&& \bigl(v^{\hat\psi}-g^{\hat\lambda} \bigr) (x,m)
\\
&&\qquad = -\int
_{\hat\psi(m)}^x \biggl(g_x(
\xi,m)-g_x \bigl(\xi,b(\xi) \bigr) -\int_{\hat\psi(m)}^\xi
\frac{g_m(y,b(y))}{
b(y)-y} b'(dy) \biggr)\,d\xi
\\
&&\qquad = \int_{\hat\psi(m)}^x\int_{\hat\psi(m)}^\xi
\biggl(g_{xm} \bigl(\xi,b(y) \bigr)+\frac{g_{m} (y,b(y) )}{b(y)-y}
\biggr)b'(dy)\,d\xi
\\
&&\qquad = \int_{\hat\psi(m)}^x \biggl(\int
_y^x g_{xm}\bigl(\xi,b(y)\bigr) +
\frac{g_m(y,b(y))}{b(y)-y} \biggr)\,d\xi\, b'(dy)
\\
&&\qquad = \int_{\hat\psi(m)}^x \bigl(b(y)-x\bigr) \biggl(
\frac{g_m(x,b(y))}{b(y)-x} -\frac{g_m(y,b(y))}{b(y)-y} \biggr)b'(dy) \ge 0,
\end{eqnarray*}
where the last inequality follows from the nondecrease of $b$ and
$x\longmapsto g_m(x,m)/\break  (m-x)$ (Assumption \ref{assB}), together with the fact
that $b(y)\ge x$ for $\hat\psi(m)\le y\le\break x\le m$.\quad\qed
\end{longlist}\noqed
\end{pf}

\begin{pf*}{Proof of Theorem~\ref{thmLookback}}
To complete the
proof of the theorem, it remains to prove that
\begin{eqnarray*}
\inf_{\lambda\in\Lambda^\mu} \bigl\{\mu(\lambda)+u^\lambda(X_0,X_0)
\bigr\} &\ge& \E^{\PP_0}_{X_0,X_0}\bigl[g \bigl(X_{\hat\tau},
X_{\hat\tau}^* \bigr)\bigr].
\end{eqnarray*}
To see this, we use the fact that the stopping time $\hat\tau$ defined
in {(\ref{tau*})} is a solution of the Skorokhod embedding
problem; that is,
$X_{\hat\tau}\sim\mu$ and $(X_{t\wedge\hat\tau})_{t\ge0}$ is a
uniformly integrable martingale; see Az\'ema and Yor
\cite{AzemaYor1,AzemaYor2}. Moreover $X^*_{\hat\tau}$ is integrable.
Then, for all $\lambda\in\Lambda^\mu$, it follows from the
definition of $u^\lambda$ that $u^\lambda(X_0,X_0)\ge
J(\lambda,\hat\tau)$, and therefore
\begin{eqnarray*}
\mu(\lambda)+u^\lambda(X_0,X_0) &\ge& \mu(
\lambda)+\E^{\PP_0}_{X_0,X_0} \bigl[g\bigl(X_{\hat\tau},
X^*_{\hat
\tau}\bigr)-\lambda(X_{\hat\tau}) \bigr]
\\
&=&
\E^{\PP_0}_{X_0,X_0} \bigl[g\bigl(X_{\hat\tau},
X^*_{\hat\tau}\bigr) \bigr].
\end{eqnarray*}\upqed
\end{pf*}

We conclude this section by a formal justification
that the function $b^{-1}$ appears naturally if one searches for the
best upper bound in {(\ref{upperboundU0})}.

\begin{longlist}[\textit{Step} 2.]
\item[\textit{Step} 1.] using expression {(\ref
{vpsi})} of
$v^\psi$, we directly compute that
\begin{eqnarray*}
&& \mu(\lambda)+u^\lambda(X_0,X_0)
\\
&&\qquad = \mu \bigl(g(
\cdot,X_0) \bigr) +\mu \bigl(g^\lambda(\cdot,X_0)
\bigr) -\int_{\psi(X_0)}^{X_0}g^\lambda_{xx}(
\xi,X_0) (X_0-\xi)\,d\xi
\\
&&\qquad = \mu \bigl(g(\cdot,X_0) \bigr) +\int g^\lambda_{xx}(
\xi,X_0) \bigl(c(\xi)-c_0(\xi){\mathbf1}_{\{\xi\le\psi(X_0)\}}
\bigr)\,d\xi
\\
&&\qquad = \mu \bigl(g(\cdot,X_0) \bigr) +\int g^\lambda_{xx}
\bigl(\xi,\psi^{-1}(\xi) \bigr) \bigl(c(\xi)-c_0(\xi){
\mathbf1}_{\{\xi\le\psi(X_0)\}} \bigr)\,d\xi
\\
&&\quad\qquad{} +\int \bigl(g_{xx}(\xi,X_0) -g_{xx} \bigl(
\xi,\psi^{-1}(\xi) \bigr) \bigr) \bigl(c(\xi)-c_0(\xi){
\mathbf1}_{\{\xi\le\psi(X_0)\}} \bigr)\,d\xi,
\end{eqnarray*}
where the second equality follows from two integrations by parts
together with the fact that $\int x\mu(dx)=X_0$; see step 1 of the
proof of Lemma 3.2 of Galichon, Henry-Labord\`ere and Touzi \cite
{ght}. Then, by using ODE {(\ref{ODE})}
satisfied by $\psi$ to change variables in the last integral, we see
that
\begin{eqnarray*}
&& \mu(\lambda)+u^\lambda(X_0,X_0)
\\
&&\qquad = \mu \bigl(g(
\cdot,X_0) \bigr) +\int \bigl\{-\gamma \bigl(\psi(m),m \bigr) +G
\bigl(\psi(m),m \bigr)\psi'(m) \bigr\} \delta \bigl(\psi(m),m
\bigr)\,dm,
\end{eqnarray*}
where we denoted
\[
\delta(x,m) := c(x)-c_0(x){\mathbf1}_{\{m\le X_0\}},
\qquad c_0(x):=(X_0-x)^+
\]
and
\[
G(x,m) := g_{xx}(x,X_0)-g_{xx}(x,m).
\]

\item[\textit{Step} 2.] The expression of $\mu(\lambda)+v^\psi$
derived in the previous step only involves the function
$\psi\in\Psi^\lambda$. Forgetting about all constraints on $\psi$,
we treat our minimization problem as a standard problem of calculus
of variations. The local Euler--Lagrange equation for this problem
is
\begin{eqnarray*}
\frac{d}{dx}(G\delta) (\psi,m) &=& -(\gamma\delta)_x(\psi,m)
+(G\delta)_x(\psi,m)\psi'.
\end{eqnarray*}
Since $(G\delta_m)(x,m)=0$, this reduces to
\begin{eqnarray*}
0 &=& (G_m\delta+\gamma\delta_x+\gamma_x
\delta) (\psi,m)
\\
&=& (m-\psi)\gamma(\psi,m) \frac{\partial}{\partial x} \biggl\{\frac{\delta(x,m)}{m-x} \biggr
\} _{x=\psi}.
\end{eqnarray*}
This shows formally that the solution of the minimization problem
\begin{eqnarray*}
&&\min_{\xi<m}\frac{\delta(x,m)}{m-x}
\end{eqnarray*}
provides a solution to the local Euler--Lagrange equation. Finally,
the solution of the above minimization problem, as recalled in
Section~\ref{secSEP}, is known to be given
by the right inverse barycenter function $b^{-1}$; see also the proof of
Lemma 3.3 in~\cite{ght}.
\end{longlist}

\begin{appendix}\label{append}
%s7 #&#
\section{Proof of Lemma \texorpdfstring{\protect\ref{27565734534563}}{3.4}}\label{sectAppendixAA}

Note that we may restrict to optimizing in {(\ref
{eqglobaloptimizationproblem})} over $\ub^{-1}(m)\leq\zeta_1\leq
\cdots\leq\zeta_n\leq m$. Indeed, fix some $\zeta$ as in {(\ref
{eqglobaloptimizationproblem})} with $\zeta_1=\cdots=\zeta_i<\min
\{\zeta_{i+1}, \ub^{-1}(m)\}$. Then all the terms featuring $\zeta
_j$ for $j\leq i$ reduce simply to $c_i(\zeta_i)/(m-\zeta_i)$ which
is nonincreasing for $\zeta_i\leq b^{-1}_i(m)$ by the discussion of
{(\ref{eqone-dimoptimisation})}. It follows that we may only
decrease the value of the objective in {(\ref
{eqglobaloptimizationproblem})} by setting $\zeta_1=\cdots=\zeta
_i=\min\{\zeta_{i+1}, \ub^{-1}(m)\}$. In consequence the problem in
{(\ref{eqglobaloptimizationproblem})} reduces to minimization
of a continuous function in a compact subset of $\mathbb{R}^n$ and
admits a minimizer $\zeta^*(m)$ with $\zeta^*_1(m)\geq\ub^{-1}(m)$.

We finally verify that $m\longmapsto\zeta^*(m)$ can be chosen to be
measurable. Indeed, for a fixed $m$, the set $F(m)$ of minimizers in
{(\ref{eqglobaloptimizationproblem})} is closed, and for any
closed $K\subset\mathbb{R}^n$, $\{m\dvtx  F(m)\cap K\neq\varnothing\}$ is
equal to $\{m\dvtx  C(m)=C_K(m)\}$ where $C_K$ is given as $C$ in {(\ref
{eqglobaloptimizationproblem})} but with a further requirement that
$(\zeta_1,\ldots,\zeta_n)\in K$. Both $C$ and $C_K$ can be obtained
through countable pointwise minimization of continuous functions and
hence are measurable, as is $\{m\dvtx  C(m)=C_K(m)\}$. Existence of a
measurable selector for $F$ now follows from Kuratowski and
Ryll-Nardzewski measurable selection theorem; see, for example, Wagner
\cite{Wagner77}, Theorem~4.1.

%s8 #&#
\section{On local martingales with simple integrands}

We now report a characterization of martingales defined as stochastic
integrals with simple integrands, which was used in Remark~\ref
{rksimplephi} and the proof of Theorem~\ref{thm-nmarginal}.

Let $\{X_t,t\in[0,T]\}$ be a $ (\mathbb{P},\{{\mathcal F}_t\}
_{0\le t\le T} )$-martingale, $0=t_0<t_1<\cdots<t_n=T$ a partition
of $[0,T]$ and $H_{t_i}$ an ${\mathcal F}_{t_i}$-measurable r.v. for
all $i=0,\ldots,n-1$. Our interest is in the process
\[
Y_t := \sum_{i=1}^n
H_{t_{i-1}} ( X_{t\wedge t_i}-X_{t\wedge t_{i-1}} ),\qquad t\in[0,T].
\]
Since $X$ is a martingale, it follows that $Y$ is a local martingale.
The following result is an easy adaptation of a similar result for
finite discrete-time local martingales reported in Jacod and Shiryaev
\cite{JacodProtter}.

%
%le8.1 #&#
\begin{Lemma} \label{lemlocalmart}
The process $Y$ is a martingale if and only if $Y_T^-$ is integrable.
\end{Lemma}

\begin{pf}
The necessary condition is obvious. We now assume $Y_T^-$ is integrable
and consider the sequence of stopping times
\begin{eqnarray*}
&&\tau_k := T\wedge\min \bigl\{t_i\dvtx  0\le i\le n, \llvert
H_{t_i}\rrvert \ge k \bigr\}, \qquad k\in\mathbb{N}.
\end{eqnarray*}
Clearly, $(\tau_k)_{k\ge1}$ is a localizing sequence for the local
martingale $Y$, taking values in the finite set $\{t_0,\ldots,t_n\}$.

We first show that $Y_t$ is integrable for all $t\in[0,T]$. By the
Jensen inequality, we have
\begin{eqnarray*}
&&Y_t^- \le \E \bigl[Y_T^-\mid {\mathcal F}_t
\bigr]\qquad\mbox{on } \{\tau_k > t_{n-1}\}.
\end{eqnarray*}
This shows that $\E [Y_t^- ]<\infty$ by sending $k\to\infty
$. We continue estimating
\begin{eqnarray*}
\E \bigl[Y_t^+ \bigr] &=& \E \Bigl[\liminf_{k\to\infty}Y_{t\wedge\tau_k}^+
\Bigr] \le \liminf_{k\to\infty} \E \bigl[Y_{t\wedge\tau_k}^+ \bigr]
\\
&=& \liminf_{k\to\infty} \E \bigl[Y_{t\wedge\tau_k}+Y_{t\wedge\tau_k}^-
\bigr]
\\
&=& \liminf_{k\to\infty}\E \bigl[Y_{t\wedge\tau_k}^- \bigr]
\\
&\le& \sum_{i=0}^n\E
\bigl[Y_{t\wedge t_i}^- \bigr] < \infty,
\end{eqnarray*}
where we used Fatou's lemma, the fact that $Y_{.\wedge\tau_k}$ is a
martingale starting from the origin and the crucial property that the
localizing sequence takes values in the finite set $\{t_i\}_{0\le i\le
n}$. %together with a conditional Jensen argument as above.
Hence $\E\llvert  Y_t\rrvert  <\infty$ for all $t\in[0,T]$.

We next show that $Y$ satisfies the martingale property. Clearly, it is
sufficient to prove the martingale property on each interval
$[t_{i-1},t_i]$. For $t_{i-1}\le s\le t\le t_i$, it follows from the
martingale property of the stopped process $Y_{.\wedge\tau_k}$,
together with the fact that the localizing sequence takes values in the
finite set $\{t_i\}_{0\le i\le n}$, that $\E[Y_t\mid{\mathcal F}_s]=Y_s$
on $\{\tau_k>s\}$. The required result follows immediately by sending
$k\to\infty$.
\end{pf}

%s9 #&#
\section{Proof of Lemma \texorpdfstring{\protect\ref{111122222454}}{5.8}}\label{sectAppendix}

We start with the computation of $\gamma_i(\psi_i,\cdot)$, as defined in
{(\ref{odemultiple})}, in terms of $\phi$ and the $\psi_i$'s.

%
%le9.1 #&#
\begin{Lemma}\label{lem-hi}
For all $i<n$, we have
$\gamma_i (\psi_i(m),m )=\frac{\phi'(m)}{
m-\psi_i(m)}
{\mathbf1}_{\{\psi_i<\overline{\psi}_{i+1}\}}$.
\end{Lemma}

\begin{pf} By direct differentiation of {(\ref
{ui})}, we see that
\begin{eqnarray*}
\partial_m v_{i-1}(x,m) &=& \partial_mv_i
\bigl(x\wedge\psi_i(m),m\bigr)
\\
&&{} +\bigl(x-\psi_i(m)\bigr)^+ \bigl[\partial_{xx}v_i
\bigl(\psi_i(m),m\bigr)\psi_i'(m) +
\partial_{xm}v_i\bigl(\psi_i(m),m\bigr)
\bigr].
\end{eqnarray*}
Using the ODE satisfied by $\psi_i$, this provides
%
%e9.1 #&#
\begin{eqnarray} \label{vim}
\partial_m v_{i-1}(x,m) &=& \partial_mv_i
\bigl(x\wedge\psi_i(m),m\bigr)\nonumber
\\
&&{}  -\frac{(x-\psi_i(m))^+}{
m-\psi_i(m)}
\partial_mv_i\bigl(x\wedge\psi_i(m),m\bigr)
\\
&=& \frac{m-x\vee\psi_i(m)}{
m-\psi_i(m)} \partial_mv_i\bigl(x\wedge
\psi_i(m),m\bigr).\nonumber
\end{eqnarray}
Since $\partial_m v_i$ is concave in $x$, we also compute by
differentiating this expression that
%
%e9.2 #&#
\begin{eqnarray}
\label{vixm} \qquad \partial_{mx} v_{i-1}(x,m) &=& {
\mathbf1}_{\{x<\psi_i(m)\}}\partial_{mx}v_i\bigl(x\wedge
\psi_i(m),m\bigr)
\nonumber\\[-8pt]\\[-8pt]\nonumber
&&{} +{\mathbf1}_{\{x>\psi_i(m)\}} \frac{-1}{m-\psi_i(m)}\partial_mv_i
\bigl(x\wedge\psi_i(m),m\bigr)\qquad\mbox{a.e.}
\nonumber
\end{eqnarray}
From the expression of $\gamma_i$, it follows from {(\ref{vim})} and
{(\ref{vixm})} that
\begin{eqnarray*}
\gamma_{i-1}(x,m) &=& {\mathbf1}_{\{x<\psi_i(m)\}}\gamma_i(x,m) =
\cdots= {\mathbf1}_{\{x<\overline{\psi}_i(m)\}}\gamma_n(x,m)
\\
&=& {\mathbf1}_{\{x<\overline{\psi}_i(m)\}}
\frac{\phi'(m)}{m-x}\qquad \mbox{a.e.}
\end{eqnarray*}\upqed
\end{pf}

\begin{pf*}{Proof of Lemma~\ref{111122222454}}
Recall that $\psi
\in\Psi^\lambda_n$ so that both $\psi_i$ and $\psi_{i}^{-1}$ are
continuous and increasing. For any integrable
function $\varphi$, the claim
%
%e9.3 #&#
\begin{eqnarray}\label{lemintigoal}
&& \int \varphi(\xi)\lambda_i''(d\xi)\nonumber
\\
&&\qquad  =
\int \Biggl( \frac{\varphi
(\psi_i(m))}{
m-\psi_i(m)} {\mathbf1}_{\{\psi_i(m)<\overline{\psi}_{i+1}(m)\}}\nonumber
\\
&&\hspace*{48pt}{} - \sum
_{j=i+1}^k \frac{\varphi(\psi_j(m))}{
m-\psi_j(m)} {\mathbf1}_{\{\psi_i(m)<\psi_j(m)=\overline{\psi}_{j}(m)\}}
\Biggr)\,d\phi(m)
\\
&&\quad\qquad{} + \int \varphi(\xi) \bigl[\partial_{xx}v_{k} \bigl(\xi,
\psi_i^{-1}(\xi) \bigr) - \partial_{xx}v_{k}
\bigl(\xi,\bigl(\psi_{i+1}^{-1} \vee\cdots\vee
\psi_{k}^{-1}\bigr) (\xi) \bigr) \bigr] \nonumber
\\
&&\hspace*{54pt}{}\times {\mathbf1}_{\{\psi_i^{-1}(\xi)
>(\psi_{i+1}^{-1}
\vee\cdots\vee
\psi_{k}^{-1})(\xi)\}} \,d\xi,
\nonumber
\end{eqnarray}
which will be proved below by induction, implies the required result
for $k=n$, and
uses the fact that $v_n=\phi$ is independent of $x$.

We start verifying {(\ref{lemintigoal})} for $k=i+1$. From the
expression of $v_i$ in
{(\ref{ui})}, we have
%
%e9.4 #&#
\begin{eqnarray}\label{uss0}
v_j &=& v^{\lambda}_{j+1}\qquad\mbox{on } \bigl\{x<
\psi_{j+1}(m)\bigr\}\quad\mbox{and}
\nonumber\\[-8pt]\\[-8pt]\nonumber
\partial_{xx}v_j&=&0\qquad
\mbox{on } \bigl\{x>\psi_{j+1}(m)\bigr\},
\end{eqnarray}
where $v^{\lambda}_{j}=v_j-\lambda_j$.

\begin{longlist}[\textit{Step} 1.]
\item[\textit{Step} 1.] To see that {(\ref{lemintigoal})} holds true with $k=i+1$, we
first decompose the integral so as to use the ODE satisfied by
$\psi_i$,
\begin{eqnarray*}
\int \varphi\lambda_i'' &=& -\int
\varphi(\xi)\partial_{xx}v_i^{\lambda} \bigl(\xi,\psi
_i^{-1}(\xi) \bigr)\,d\xi +\int\varphi(\xi)
\partial_{xx}v_i \bigl(\xi,\psi_i^{-1}(
\xi ) \bigr)\,d\xi
\\
&=& \int\varphi \bigl(\psi_i(m) \bigr)\gamma_i \bigl(
\psi _i(m),m \bigr)\,dm +\int\varphi(\xi)\partial_{xx}v_i
\bigl(\xi,\psi_i^{-1}(\xi ) \bigr)\,d\xi.
\end{eqnarray*}
Substitute the expression of $\gamma_i(\psi_i,\cdot)$ from Lemma
\ref{lem-hi}, and use {(\ref{uss0})} for the second integral,
\begin{eqnarray*}
\int \varphi\lambda_i'' &=& \int
\frac{\varphi(\psi_i(m))}{m-\psi_i(m)} {\mathbf1}_{\{\psi_i(m)<\overline{\psi}_{i+1}(m)\}}\,d\phi(m)
\\
&&{} +\int\varphi(\xi)
\partial_{xx}v_{i+1}^\lambda \bigl(\xi,
\psi_i^{-1}(\xi) \bigr) {\mathbf1}_{\{\psi_{i+1}^{-1}(\xi)<\psi_i^{-1}(\xi)\}}\,d\xi
\\
&=& \int\frac{\varphi(\psi_i(m))}{m-\psi_i(m)} {\mathbf1}_{\{\psi_i(m)<\overline{\psi}_{i+1}(m)\}}\,d\phi(m)
\\
&&{}+\int\varphi(\xi)
\partial_{xx}v_{i+1}^\lambda \bigl(\xi,
\psi_{i+1}^{-1}(\xi) \bigr) {\mathbf1}_{\{\psi_{i+1}^{-1}(\xi)<\psi_i^{-1}(\xi)\}}\,d\xi
\\
&& {}+\int\varphi(\xi) \bigl[\partial_{xx}v_{i+1}^\lambda
\bigl(\xi,\psi_i^{-1}(\xi) \bigr) -\partial_{xx}v_{i+1}^\lambda
\bigl(\xi,\psi_{i+1}^{-1}(\xi) \bigr) \bigr]
\\
&&\hspace*{21pt}{}\times {
\mathbf1}_{\{\psi_{i+1}^{-1}(\xi)<\psi_i^{-1}(\xi)\}}\,d\xi.
\end{eqnarray*}
Then, by using again ODE {(\ref{ODE})} satisfied by $\psi_{i+1}$
together with the expression of $\gamma_{i+1}(\psi_{i+1},\cdot)$ from
Lemma~\ref{lem-hi}, we get
\begin{eqnarray*}
\int \varphi\lambda_i'' &=& \int
\frac{\varphi (\psi_i(m) )}{m-\psi_i(m)} {\mathbf1}_{\{\psi_i(m)<\overline{\psi}_{i+1}(m)\}}\,d\phi(m)
\\
&&{} -\int\frac{\varphi (\psi_{i+1}(m) )}{m-\psi_{i+1}(m)} {
\mathbf1}_{\{\psi_i(m)<\psi_{i+1}(m)=\overline{\psi}_{i+1}(m)\}}\,d\phi(m)
\\
&&{} +\int\varphi(\xi) \bigl[\partial_{xx}v_{i+1}^\lambda
\bigl(\xi,\psi_i^{-1}(\xi) \bigr) -\partial_{xx}v_{i+1}^\lambda
\bigl(\xi,\psi_{i+1}^{-1}(\xi) \bigr) \bigr]
\\
&&\hspace*{21pt}{}\times {
\mathbf1}_{\{\psi_{i+1}^{-1}(\xi)<\psi_i^{-1}(\xi)\}}\,d\xi,
\end{eqnarray*}
which we recognize to be the required equality {(\ref{lemintigoal})}
for $k=i+1$.

\item[\textit{Step} 2.] We next assume that {(\ref{lemintigoal})} holds for some $k<n-1$,
and verify it for $k+1$. For simplicity, we denote
$\psi_{{i+1},j}^{-1}:=\psi_{i+1}^{-1}\vee\cdots\vee\psi_{j}^{-1}$.
By {(\ref{uss0})}, we compute that
\begin{eqnarray*}
A &:=& \int\varphi(\xi) \bigl[\partial_{xx}v_k \bigl(
\xi,\psi_i^{-1}(\xi) \bigr) -\partial_{xx}v_k
\bigl(\xi,\psi_{i+1,k}^{-1}(\xi) \bigr) \bigr]{
\mathbf1}_{\{\psi_i^{-1}(\xi)>\psi_{i+1,k}^{-1}(\xi)\}} \,d\xi
\\
&=& \int\varphi(\xi){\mathbf1}_{\{\psi_i^{-1}(\xi)>\psi_{i+1,k}^{-1}(\xi
)\}}
\\
&&\hspace*{8pt}{}\times \bigl[ \bigl\{\partial_{xx}v_{k+1}
\bigl(\xi,\psi_i^{-1}(\xi) \bigr) -\lambda_{k+1}''(
\xi) \bigr\}{\mathbf1}_{\{\psi_{k+1}^{-1}(\xi)<\psi_i^{-1}(\xi)\}}
\\
&&\hspace*{23pt}{}  - \bigl\{\partial_{xx}v_{k+1} \bigl(\xi,
\psi_{i+1,k}^{-1}(\xi) \bigr) -\lambda_{k+1}''(
\xi) \bigr\}{\mathbf1}_{\{\psi_{k+1}^{-1}(\xi)<\psi_{i+1,k}^{-1}(\xi)\}} \bigr] \,d\xi
\\
&=& \int\varphi(\xi){\mathbf1}_{\{\psi_i^{-1}(\xi)>\psi_{i+1,k}^{-1}(\xi
)\}}
\\
&&\hspace*{8pt}{}\times \bigl[{\mathbf1}_{\{\psi_{i+1,k}^{-1}(\xi)
<\psi_{k+1}^{-1}(\xi)
<\psi_i^{-1}(\xi)\}}
\partial_{xx}v^\lambda_{k+1} \bigl(\xi,
\psi_i^{-1}(\xi) \bigr)
\\
&&\hspace*{23pt}{}  +{\mathbf1}_{\{\psi_{k+1}^{-1}(\xi)<\psi_{i+1,k}^{-1}(\xi)\}}
\\
&&\hspace*{33pt}{}\times \bigl\{\partial_{xx}v_{k+1}
\bigl(\xi,\psi_i^{-1}(\xi) \bigr) -\partial_{xx}v_{k+1}
\bigl(\xi,\psi_{i+1,k}^{-1}(\xi) \bigr) \bigr\} \bigr] \,d\xi
\\
&=& \int\varphi(\xi){\mathbf1}_{\{\psi_i^{-1}(\xi)>\psi_{i+1,k}^{-1}(\xi
)\}}
\\
&&\hspace*{8pt}{}\times \bigl[{\mathbf1}_{\{\psi_{i+1,k}^{-1}(\xi)
<\psi_{k+1}^{-1}(\xi)
<\psi_i^{-1}(\xi)\}}
\partial_{xx}v^\lambda_{k+1} \bigl(\xi,
\psi_{k+1}^{-1}(\xi) \bigr)
\\
&&\hspace*{23pt}{}  +{\mathbf1}_{\{\psi_{i+1,k}^{-1}(\xi)
<\psi_{k+1}^{-1}(\xi)<\psi_i^{-1}(\xi)\}}
\\
&&\hspace*{33pt}{}\times \bigl\{\partial_{xx}v^\lambda_{k+1}
\bigl(\xi,\psi_i^{-1}(\xi ) \bigr) -\partial_{xx}v^\lambda_{k+1}
\bigl(\xi,\psi_{k+1}^{-1}(\xi) \bigr) \bigr\}
\\
&&\hspace*{23pt}{}  +{\mathbf1}_{\{\psi_{k+1}^{-1}(\xi)<\psi_{i+1,k}^{-1}(\xi)\}}
\\
&&\hspace*{33pt}{}\times  \bigl\{\partial_{xx}v_{k+1}
\bigl(\xi,\psi_i^{-1}(\xi) \bigr) -\partial_{xx}v_{k+1}
\bigl(\xi,\psi_{i+1,k}^{-1}(\xi) \bigr) \bigr\} \bigr] \,d\xi.
\end{eqnarray*}
Putting together the two last terms, we see that
\begin{eqnarray*}
A &=& \int\varphi(\xi){\mathbf1}_{\{\psi_i^{-1}(\xi)>\psi_{i+1,k}^{-1}(\xi
)\}}
\\
&&\hspace*{8pt}{}\times \bigl[{\mathbf1}_{\{\psi_{i+1,k}^{-1}(\xi)
<\psi_{k+1}^{-1}(\xi)
<\psi_i^{-1}(\xi)\}}
\partial_{xx}v^\lambda_{k+1} \bigl(\xi,
\psi_{k+1}^{-1}(\xi) \bigr)
\\
&&\hspace*{23pt}{} +{\mathbf1}_{\{\psi_{k+1}^{-1}(\xi)
<\psi_i^{-1}(\xi)\}}
\\
&&\hspace*{33pt}{}\times \bigl\{\partial_{xx}v^\lambda_{k+1}
\bigl(\xi,\psi_i^{-1}(\xi) \bigr) -\partial_{xx}v^\lambda_{k+1}
\bigl(\xi,\psi_{i+1,k+1}^{-1}(\xi) \bigr) \bigr\} \bigr] \,d\xi.
\end{eqnarray*}
Finally, using ODE {(\ref{odemultiple})} satisfied by $\psi
_{k+1}$ in the
first term, together with the expression of
$\gamma_{k+1}(\psi_{k+1},\cdot)$ from Lemma~\ref{lem-hi}, we see that
\begin{eqnarray*}
A &=& -\int\varphi \bigl(\psi_{k+1}(m) \bigr) \frac{\varphi (\psi_{k+1}(m) )}{
\psi_{k+1}(m)-m} {
\mathbf1}_{\{\psi_i(m)<\psi_{k+1}(m)=\overline{\psi}_{k+1}(m)\}} \,d\phi(m)
\\
&&{} +\int\varphi(\xi) \bigl[\partial_{xx}v^\lambda_{k+2}
\bigl(\xi,\psi_i^{-1}(\xi) \bigr) -\partial_{xx}v^\lambda_{k+2}
\bigl(\xi,\psi_{i+1,k+2}^{-1}(\xi) \bigr) \bigr]
\\
&&\hspace*{20pt}{}\times {
\mathbf1}_{\{\psi_i^{-1}(\xi)>\psi_{i+1,k+1}^{-1}(\xi)\}} \,d\xi,
\end{eqnarray*}
which is precisely the required expression in order to justify that
{(\ref{lemintigoal})} holds for $k+1$.\quad\qed
\end{longlist}\noqed
\end{pf*}

%s10 #&#
\section{Proofs for statements in Section~\texorpdfstring{\protect\ref{55555222222222}}{4.2}}
\label{appproofsiAYembedding}

\begin{pf*}{Proof of Lemma~\ref{propPathwiseEquality}}
Fix $m > Z_0$, and write for notational convenience $\zeta_i = \zeta
_i(m),  \eta_i = \eta_i(m)$.
Let $n_1<\cdots<n_k = n$ be such that
\[
\zeta_1 = \cdots= \zeta_{n_1} < \zeta_{n_1 + 1} =
\cdots= \zeta _{n_2} < \cdots< \zeta_{n_{k-1}+1} = \cdots=
\zeta_{n_k} = \zeta_n = \eta_n.
\]
Then by
\[
Z_{t_{j}} \geq\zeta_{j} \quad\Longrightarrow\quad
Z_{t_l} \geq \zeta_{j} \qquad\forall l \geq j
\]
for all $j \leq n$, which follows directly from the definitions {\rm
(\ref{eqiteratedAY})} and {(\ref{eqchangeoftimeAY})},
we obtain by identifying the terms as in the proof of Proposition~\ref
{lemTrajectorialInequalityOrderedCase},
\begin{eqnarray*}
\Upsilon_n(Z, m, \zeta) &= &\sum_{j=1}^{k}
\biggl( \frac{ (Z_{t_{n_j}}-\zeta_i
)^+}{m-\zeta_{n_j}} +{\mathbf1}_{ \{ {Z^*_{t_{n_{j-1}}} < m \leq Z^*_{t_{n_j}}}  \}
} \frac{ m - Z_{t_{n_j}}}{m-\zeta_{n_j}} \biggr)
\\
&&{} -\sum_{j=1}^{k-1} \biggl(
\frac{ (Z_{t_{n_j}}-\zeta
_{n_{j+1}} )^+}{m-\zeta_{n_{j+1}}} +{\mathbf1}_{ \{ {m \leq Z^*_{t_{n_j}}, \zeta_{n_{j+1}} \leq
Z_{t_{n_j}} }  \} } \frac{Z_{t_{n_{j+1}}} - Z_{t_{n_j}}}{m-\zeta_{n_{j+1}}} \biggr).
\end{eqnarray*}
Therefore,
%in view of \citet[Lemma 3.1]{OblSp13},
it is enough to prove the claim for
%\begin{eqnarray*}
%\jmath_i = i-1 \qquad\forall i=1,\ldots,n,
%\end{eqnarray*}
%i.e. for
the case
\[
\zeta_1 = \eta_1 <\zeta_2 =
\eta_2 <\cdots<\zeta_n = \eta_n.
\]

By the same induction as in the proof of Proposition~\ref
{lemTrajectorialInequalityOrderedCase} it remains to prove that
$(Z,Z^*)$ achieves equality in {(\ref{eqpwproofindstart})},
{(\ref{eqrebalancingcalendarspreads})} and {(\ref
{eqforwardtransaction})}.
As for equality in {(\ref{eqpwproofindstart})} we note that
\[
Z^*_{t_1} \geq m \quad\Longrightarrow\quad Z_{t_1} \geq
\zeta_1 \quad\mbox{and}\quad Z^*_{t_1} < m \quad
\Longrightarrow\quad Z_{t_1} \leq\zeta_1.
\]
Equality in {(\ref{eqrebalancingcalendarspreads})} holds by
\begin{eqnarray*}
Z^*_{t_{n-1}} \geq m,\qquad Z_{t_{n-1}} \geq\zeta_n &\quad\Longrightarrow\quad& Z_{t_n} \geq\zeta_n,
\\
Z^*_{t_{n-1}} \geq m,\qquad Z_{t_{n-1}} < \zeta_n &\quad\Longrightarrow \quad& Z_{t_n} < \zeta_n,
\end{eqnarray*}
which one verifies using the definition of the iterated Az\'
ema--Yor-type embedding.
Similarly, equality in {(\ref{eqforwardtransaction})} holds by
\begin{eqnarray*}
Z^*_{t_{n-1}} < m,\qquad Z^*_{t_n} \geq m &\quad\Longrightarrow\quad&
Z_{t_n} \geq\zeta_n,
\\
Z^*_{t_{n-1}} < m,\qquad Z^*_{t_n} < m &\quad\Longrightarrow\quad&
Z_{t_n} \leq\zeta_n.
\end{eqnarray*}
The claim follows.
\end{pf*}

\begin{pf*}{Proof of equation {(\ref{equniqueminimizers})}}
Finally, we argue that {(\ref{equniqueminimizers})} holds under
Assumption $\circledast$ of Ob\l\'oj and Spoida \cite{OblSp}.
Let $\tilde{\zeta}_i(m):= \min_{j\geq i} \eta_j(m)$.
%By definition of the \iAY type embedding, see Section~\ref{55555222222222},
%we have
%\begin{eqnarray*}
%Z_{t_j} = Z_{t_{j+1}} = \dots= Z_{t_n} \quad\mbox{and} \quad
%Z^*_{t_j} = Z^*_{t_{j+1}} = \dots= Z^*_{t_n} \qquad\mbox{on the set $
%\left\{ Z_{t_j} = \tilde{\zeta}_j(\barr{Z}_{t_j}) \right\}$.}
%\end{eqnarray*}

First fix $m \in(X_0, r^{\mu_n})$ such that $\tilde{\zeta}_n(m)<m$.
Then Lemma~\ref{propPathwiseEquality} and Proposition~\ref
{propmartineq} yield
\[
\PP \bigl[ Z^*_{t_n} \geq m \bigr] = \E \bigl[ \Upsilon_n(Z,m,
\tilde{\zeta}) \bigr] \leq C(m).
\]
Next we show that under Assumption $\circledast$ this
inequality is strict whenever $\zeta^*(m) \neq\tilde{\zeta}(m)$.
This will be a contradiction to the optimality of $\zeta^*$ since, by
invoking the embedding property $Z_{t_{i}} \sim\mu_i$, we must have
$\E [ \Upsilon_n(Z,m,\tilde{\zeta})  ] \geq C(m)$.

\begin{longlist}[\textit{Case} A.]
\item[\textit{Case }A.]   $\tilde{\zeta}_j(m) > \zeta^*_j(m)$.
Assume initially that $m \in(X_0, r^{\mu_j}]$.
Then on the set
\[
\bigl\{ Z_{t_{j}} > \zeta^*_j(m), Z^*_{t_j} < m
\bigr\} \supseteq \bigl\{ Z^*_{t_{j}} \in\mathcal{O}_A ,
Z^*_{t_n} \in\mathcal{O}_A \bigr\} =: \mathcal{Z}_A,
\]
where $\mathcal{O}_A \subseteq[X_0,m)$ is a (suitable) open interval,
we obtain
\[
\Upsilon_n\bigl( Z, m, \zeta^* \bigr) \geq\Upsilon_j
\bigl( Z, m, \zeta^* \bigr) \stackrel{\fontsize{8.36pt}{10pt}{{(\ref{eqforwardtransaction})}}}{>} {
\mathbf1}_{ \{ {m \leq Z^*_{t_j} }  \} } = {\mathbf1}_{ \{
{m \leq Z^*_{t_n} }  \} } \stackrel{\fontsize{8.36pt}{10pt}{\mathrm{Lemma}~\ref{propPathwiseEquality}}} {=} \Upsilon_n( Z, m, \tilde{\zeta}).
\]
Note that $\E [ \Upsilon_n( Z, m, \zeta^* )  ] = C(m)$.
It follows that if $\PP [ \mathcal{Z}_A  ] > 0 $, then\break $\E
 [ \Upsilon_n( Z, m, \tilde{\zeta})  ] < C(m)$ as
required. However, this is clear by the assumption that $m \in(X_0,
r^{\mu_j}]$ and elementary properties of Brownian motion.

If $m > r^{\mu_j}$, then by Lemma~\ref{propPathwiseEquality} we get
$\E [ \Upsilon_j( Z, m, \tilde{\zeta})  ] = 0$. If $\sum_{i=1}^j  (\frac{c_i(\zeta^*_i)}{m-\zeta^*_i} -\frac{c_i(\zeta
^*_{i+1})}{m- \zeta^*_{i+1}} {\mathbf1}_{\{i<n\}}  ) > 0$, then $\E
 [\Upsilon_j( Z, m, \zeta^*)  ] > 0$ and hence
\begin{eqnarray*}
\E \bigl[ {\mathbf1}_{ \{ {Z_{t_n} < m }  \} } \Upsilon_n\bigl( Z, m, \zeta^*
\bigr) \bigr] &\geq& \E \bigl[ {\mathbf1}_{ \{ {Z_{t_n} < m }
 \} } \Upsilon_j\bigl(
Z, m, \zeta^*\bigr) \bigr] > 0
\\
&=& \E \bigl[ {\mathbf1}_{ \{ {Z_{t_n} < m }  \} }
\Upsilon_n( Z, m, \tilde{\zeta}) \bigr].
\end{eqnarray*}
If $\sum_{i=1}^j  (\frac{c_i(\zeta^*_i)}{m-\zeta^*_i} -\frac
{c_i(\zeta^*_{i+1})}{m- \zeta^*_{i+1}} {\mathbf1}_{\{i<n\}}  ) =
0$, then a contradiction to Assumption $\circledast$ is
obtained if $\zeta^*_j(m) < \tilde{\zeta}_j(m)$.

\item[\textit{Case} B.] $\tilde{\zeta}_j(m) < \zeta^*_j(m)$.
We can, without loss of generality, take $m \in(X_0, r^{\mu_j}]$.
Indeed, if this is not the case, then we set $j' \geq j$ such that
$\tilde{\zeta}_j(m) = \tilde{\zeta}_{j+1}(m) = \cdots= \tilde
{\zeta}_{j'}(m) = \eta_{j'}(m) < \tilde{\zeta}_{j'+1}(m)$.
For this $j'$ we then have $\tilde{\zeta}_{j'}(m) < \zeta^*_{j'}(m)$
and $m \in(X_0, r^{\mu_{j'}}]$ as $\eta_{j'}(m)<m$.

Then on the set
\[
\bigl\{ Z_{t_j} < \zeta^*_j(m), Z^*_{t_j} \geq
m \bigr\} \supseteq \bigl\{ Z^*_{t_j} \in\mathcal{O}_B
\bigr\} =: \mathcal{Z}_B,
\]
where $\mathcal{O}_B \subseteq[m,\infty)$ is a (suitable) open interval,
we obtain
\[
\Upsilon_n\bigl( Z, m, \zeta^* \bigr) \geq\Upsilon_j
\bigl( Z, m, \zeta^* \bigr) \stackrel{\fontsize{8.36pt}{10pt}{(\ref{eqforwardtransaction})}} {>} 1 = {\mathbf
1}_{ \{ {m \leq Z^*_{t_j} }  \} } = {\mathbf1}_{ \{ {m
\leq Z^*_{t_n} }  \} } \stackrel{\fontsize{8.36pt}{10pt}{\mathrm{Lemma}~\ref{propPathwiseEquality}}} {=} \Upsilon_n( Z, m, \tilde{\zeta} ),
\]
which yields a contradiction in a similar fashion to case A because
$\PP [ \mathcal{Z}_B  ]>0$.

Now consider $\tilde{\zeta}_n(m) = m$.
We assume for simplicity of the argument that $\tilde{\zeta
}_{n-1}(m)<\tilde{\zeta}_{n}(m)$. (By, e.g., use of Lemma A.1 of Ob\l
\'oj and Spoida \cite{OblSp}, the general case can be reduced to this case.)
Then optimality of $\zeta^*_n(m)$ yields
\[
\frac{ c_n(z) - c_{n-1}(z) }{ m-
z } \bigg| _{z=\zeta^*_n(m)} \leq \frac{ c_n(z) - c_{n-1}(z) }{ m-
z } \bigg|
_{z = \tilde{\zeta}_n(m) = m } .
\]
From this a direct contradiction to Assumption $\circledast$
is obtained if $\zeta^*_n(m) \neq\tilde{\zeta}_n(m)$.
Hence $\zeta^*_n(m) = \tilde{\zeta}_n(m)$.
Then $\zeta^*_1,\ldots,\zeta^*_{n-1}$ are also optimal for the $n-1$
marginal problem {(\ref{eqglobaloptimizationproblem})}.
Further, the $n-1$ marginal Az\'ema--Yor embedding coincides with the
$n$ marginal iterated Az\'ema--Yor embedding until time $t_{n-1}$.
Hence, by induction, $\zeta^*_i(m) = \tilde{\zeta}_i(m)$ for all $i
< n$.\quad\qed
\end{longlist}\noqed
\end{pf*}
\end{appendix}

\section*{Acknowledgments}
The title is borrowed from Hobson \cite{Hobsonmax} and Brown, Hobson
and Rogers \cite{brownhobsonrogers}.
%Jan Ob\l\'oj thankfully acknowledges support from the ERC Starting
%Grant {\sc RobustFinMath} 335421, the Oxford-Man Institute of
%Quantitative Finance and St John's College in Oxford. Peter Spoida
%gratefully acknowledges scholarships from the Oxford-Man Institute of
%Quantitative Finance and the DAAD. Nizar Touzi gratefully acknowledges
%financial support from the ERC Advanced Grant 321111 ROFIRM, the Chair
%{\it Financial Risks} of the {\it Risk Foundation} sponsored by Soci\'
%et\'e G\'en\'erale and the Chair {\it Finance and Sustainable
%Development} sponsored by EDF and CA-CIB.

%\begin{supplement}[id=suppA]
%\slink[doi]{10.1214/00-AAPXXXXSUPP} %[doi,text={...}] - jei reikia
%suskaldyti doi
%\sfilename{aapXXXX\_supp.pdf}
%\end{supplement}

% imsref loaded by linak, 2014-12-22 12:19:51
% imsref loaded by linak, 2015-02-03 09:48:54
% imsref loaded by linak, 2015-02-03 09:53:45

\printaddresses
\end{document}